\documentclass{article}
\usepackage{amssymb}
\title{{\bf A logarithmic generalization of tensor product theory
 for modules for a vertex operator algebra}}
\author{Yi-Zhi Huang, James Lepowsky and Lin Zhang}
\date{}
\begin{document}
\bibliographystyle{alpha}
\maketitle

\newtheorem{rema}{Remark}[section]
\newtheorem{propo}[rema]{Proposition}
\newtheorem{theo}[rema]{Theorem}
\newtheorem{defi}[rema]{Definition}
\newtheorem{lemma}[rema]{Lemma}
\newtheorem{corol}[rema]{Corollary}
\newtheorem{exam}[rema]{Example}
\newtheorem{nota}[rema]{Notation}
\newcommand{\ba}{\begin{array}}
\newcommand{\ea}{\end{array}}
\newcommand{\be}{\begin{equation}}
\newcommand{\ee}{\end{equation}}
\newcommand{\bea}{\begin{eqnarray}}
\newcommand{\eea}{\end{eqnarray}}
\newcommand{\nno}{\nonumber}
\newcommand{\lbar}{\bigg\vert}
\newcommand{\p}{\partial}
\newcommand{\dps}{\displaystyle}
\newcommand{\bra}{\langle}
\newcommand{\ket}{\rangle}
\newcommand{\res}{\mbox{\rm Res}}
\newcommand{\wt}{\mbox{\rm wt}\;}
\newcommand{\pf}{{\it Proof}\hspace{2ex}}
\newcommand{\epf}{\hspace{2em}$\square$}
\newcommand{\epfv}{\hspace{1em}$\square$\vspace{1em}}
\newcommand{\ob}{{\rm ob}\,}
\renewcommand{\hom}{{\rm Hom}}

\newcommand{\dlt}[3]{#1 ^{-1}\delta \bigg( \frac{#2 #3 }{#1 }\bigg) }
\newcommand{\dlti}[3]{#1 \delta \bigg( \frac{#2 #3 }{#1 ^{-1}}\bigg) }

\makeatletter
\newlength{\@pxlwd} \newlength{\@rulewd} \newlength{\@pxlht}
\catcode`.=\active \catcode`B=\active \catcode`:=\active \catcode`|=\active
\def\sprite#1(#2,#3)[#4,#5]{
   \edef\@sprbox{\expandafter\@cdr\string#1\@nil @box}
   \expandafter\newsavebox\csname\@sprbox\endcsname
   \edef#1{\expandafter\usebox\csname\@sprbox\endcsname}
   \expandafter\setbox\csname\@sprbox\endcsname =\hbox\bgroup
   \vbox\bgroup
  \catcode`.=\active\catcode`B=\active\catcode`:=\active\catcode`|=\active
      \@pxlwd=#4 \divide\@pxlwd by #3 \@rulewd=\@pxlwd
      \@pxlht=#5 \divide\@pxlht by #2
      \def .{\hskip \@pxlwd \ignorespaces}
      \def B{\@ifnextchar B{\advance\@rulewd by \@pxlwd}{\vrule
         height \@pxlht width \@rulewd depth 0 pt \@rulewd=\@pxlwd}}
      \def :{\hbox\bgroup\vrule height \@pxlht width 0pt depth
0pt\ignorespaces}
      \def |{\vrule height \@pxlht width 0pt depth 0pt\egroup
         \prevdepth= -1000 pt}
   }
\def\endsprite{\egroup\egroup}
\catcode`.=12 \catcode`B=11 \catcode`:=12 \catcode`|=12\relax
\makeatother

\def\hboxtr{\FormOfHboxtr} 
\sprite{\FormOfHboxtr}(25,25)[0.5 em, 1.2 ex] 

:BBBBBBBBBBBBBBBBBBBBBBBBB |
:BB......................B |
:B.B.....................B |
:B..B....................B |
:B...B...................B |
:B....B..................B |
:B.....B.................B |
:B......B................B |
:B.......B...............B |
:B........B..............B |
:B.........B.............B |
:B..........B............B |
:B...........B...........B |
:B............B..........B |
:B.............B.........B |
:B..............B........B |
:B...............B.......B |
:B................B......B |
:B.................B.....B |
:B..................B....B |
:B...................B...B |
:B....................B..B |
:B.....................B.B |
:B......................BB |
:BBBBBBBBBBBBBBBBBBBBBBBBB |

\endsprite

\def\shboxtr{\FormOfShboxtr} 
\sprite{\FormOfShboxtr}(25,25)[0.3 em, 0.72 ex] 

:BBBBBBBBBBBBBBBBBBBBBBBBB |
:BB......................B |
:B.B.....................B |
:B..B....................B |
:B...B...................B |
:B....B..................B |
:B.....B.................B |
:B......B................B |
:B.......B...............B |
:B........B..............B |
:B.........B.............B |
:B..........B............B |
:B...........B...........B |
:B............B..........B |
:B.............B.........B |
:B..............B........B |
:B...............B.......B |
:B................B......B |
:B.................B.....B |
:B..................B....B |
:B...................B...B |
:B....................B..B |
:B.....................B.B |
:B......................BB |
:BBBBBBBBBBBBBBBBBBBBBBBBB |

\endsprite

\begin{abstract}

We describe a logarithmic tensor product theory for certain module
categories for a ``conformal vertex algebra.'' In this theory, which
is a natural, although intricate, generalization of earlier work of
Huang and Lepowsky, we do not require the module categories to be
semisimple, and we accommodate modules with generalized weight
spaces. The corresponding intertwining operators contain logarithms of
the variables.
\end{abstract}

\renewcommand{\theequation}{\thesection.\arabic{equation}}
\renewcommand{\therema}{\thesection.\arabic{rema}}
\setcounter{equation}{0}
\setcounter{rema}{0}

\section{Introduction}

In this paper, we present a generalization of the tensor product
theory developed by the first two authors in
\cite{tensorAnnounce}--\cite{tensor3}, \cite{tensor4} and
\cite{diff-eqn} to a logarithmic tensor product theory for certain
module categories for a ``conformal vertex algebra.'' We focus on
explaining new results and sophisticated techniques in this
logarithmic theory, rather than on the full details of the
proofs, which are being presented
in a longer paper \cite{HLZ}.

Logarithmic structure in conformal field theory was first introduced
by physicists to describe disorder phenomena \cite{Gu}. A lot of
progress has been made by physicists on this subject. We refer the
interested reader to review articles \cite{Ga}, \cite{F} and \cite{RT}
(and references therein). Such structures also arise naturally in the
representation theory of vertex operator algebras. In fact, in the
construction of intertwining operator algebras, the first author
proved (see \cite{diff-eqn}) that if modules for the vertex operator
algebra satisfy a certain cofiniteness condition, then products of the
usual intertwining operators satisfy certain systems of differential
equations with regular singular points. In addition, it was proved in
\cite{diff-eqn} that if the vertex operator algebra satisfies certain
finite reductivity conditions, then the analytic extensions of
products of the usual intertwining operators have no logarithmic
terms.

In the case when the vertex operator algebra satisfies only the
cofiniteness condition but not the finite reductivity conditions, the
products of intertwining operators still satisfy systems of
differential equations with regular singular points.  But in this
case, the analytic extensions of such products of intertwining
operators might have logarithmic terms. This means that if we want to
generalize the results in \cite{tensorAnnounce}--\cite{tensor3},
\cite{tensor4} and \cite{diff-eqn} to the case in which the finite
reductivity properties are not always satisfied, we have to consider
intertwining operators involving logarithmic terms.

In \cite{Mi}, Milas introduced and studied what he called
``logarithmic modules'' and ``logarithmic intertwining operators.''
Roughly speaking, logarithmic modules are weak modules for vertex
operator algebras that are direct sums of generalized eigenspaces for
the operator $L(0)$ (we will start to call such weak modules
``generalized modules'' in this paper), and logarithmic intertwining
operators are operators which depend not only on powers of a (formal
or complex) variable $x$, but also on its logarithm $\log x$.

In \cite{HLZ}, we have generalized the tensor product theory developed
by the first two authors to the category of generalized modules for a
``conformal vertex algebra'', or even more generally, for a ``M\"obius
vertex algebra'', satisfying suitable conditions.  The special
features of the logarithm function make the logarithmic theory very
subtle and interesting.  Although we have shown that all the main
theorems in the original tensor product theory developed by the first
two authors still hold in the logarithmic theory, many of the proofs
and techniques being used involve certain new and sophisticated
techniques and have surprising connections with certain combinatorial
identities.

In fact, one of the main technical tools that we use in this paper
(and indeed, must use) is the same as a certain critical and subtle
ingredient of the original tensor product theory of Huang and
Lepowsky: a condition called the ``compatibility condition'' in
\cite{tensorAnnounce}--\cite{tensor3} and \cite{tensor4}. In the
construction of the associativity isomorphism in \cite{tensor4}, one
of the important steps is to identify the ``intermediate module,''
which is the tensor product module of two modules.  The compatibility
condition is the most important condition among those conditions
characterizing the tensor product module and thus it is natural to use
it to identify the intermediate module.  (In fact, in a recent paper
\cite{HLLZ} we have proved by means of suitable counterexamples that
the compatibility condition is required in the theory.)  The method
based on the compatibility condition depends on the simple but
powerful result that any ``intertwining map'' determines an
intertwining operator uniquely.  To generalize the original theory to
the logarithmic theory, one first has to generalize this result, and
the generalization of this basic result involves some subtle
properties of the logarithmic function.  In particular, we see that
the method based on the compatibility condition, and therefore all the
proofs in the logarithmic theory, are more delicate than in the
original theory. In the present paper, we shall explain the results
obtained in this new logarithmic theory.  Further details and proofs
are contained in \cite{HLZ}.

One important application of this generalization is to the category
${\cal O}_\kappa$ of certain modules for an affine Lie algebra studied
by Kazhdan and Lusztig in their series of papers
\cite{KL1}--\cite{KL5}. It has been shown in \cite{Z} by the third
author that, among other things, all the conditions needed to apply
our theory to this module category are satisfied.  As a result, it is
proved in \cite{Z} that there is a natural vertex tensor category
structure on this module category, giving a new construction, in the
context of general vertex operator algebra theory, of the braided
tensor category structure on ${\cal O}_\kappa$.  The methods used in
\cite{KL1}--\cite{KL5} were very different.

The present work is a natural (although subtle) extension of
\cite{tensorAnnounce}--\cite{tensor3}, \cite{tensor4} and
\cite{diff-eqn}.  This work was presented by one of us (L.Z.) in
seminars at Rutgers and at Stony Brook in April, 2003.

\paragraph{Acknowledgment}
The authors would like to thank Sasha Kirillov Jr., Masahiko Miyamoto,
Kiyokazu Nagatomo and Bin Zhang for their early interest in our work
on logarithmic tensor product theory, and Bin Zhang in particular for
inviting L.Z. to lecture on this ongoing work at Stony Brook.
Y.-Z.H. would also like to thank Paul Fendley for raising his interest
in logarithmic conformal field theory.  The authors gratefully
acknowledge partial support from NSF grants DMS-0070800 and DMS-0401302.

\setcounter{equation}{0}
\setcounter{rema}{0}

\section{The setting}
In this section we define some notions that will be used in this
paper. We shall assume the reader has some familiarity with the
material in \cite{B}, \cite{FLM} and \cite{FHL}, such as the formal
delta function, the binomial expansion convention, the notion of
vertex (operator) algebra, and the notion of module for a vertex
(operator) algebra. We also refer the reader to \cite{LL} for such
material.

Throughout this paper we use the following usual conventional
notation: the symbols $x$, $x_0$, $x_1$, $x_2$, $\dots$ , $y$, $y_0$,
$y_1$, $y_2$, $\dots$ will denote commuting independent formal
variables, and by contrast, the symbols $z$, $z_0$, $z_1$, $z_2$,
$\dots$ will denote complex numbers in specified domains.  We use the
notation $\mathbb{N}$ for the nonnegative integers and $\mathbb{Z}_+$
for the positive integers.

We will use the following version of the notion of ``conformal vertex
algebra'': A {\em conformal vertex algebra} is a vertex algebra
equipped with a $\mathbb{Z}$-grading and with a conformal vector
$\omega$ satisfying the usual conditions. Hence the only
difference between a conformal vertex algebra and a vertex operator
algebra is that a vertex operator algebra further satisfies two
``grading restriction conditions'' (see \cite{FLM} and \cite{FHL}).

\begin{rema}{\rm
Analogous to the notion of ``quasi-vertex operator algebra'' in
\cite{FHL}, a slightly more general notion of ``M\"obius vertex
algebra'' can be defined: A {\it M\"obius vertex algebra} is a vertex
algebra equipped with a representation of
${\mathfrak{s}}{\mathfrak{l}}(2)$ on $V$ satisfying the same
conditions as those for the operators $L(-1)$, $L(0)$ and $L(1)$ in
the definition of vertex operator algebra.  Essentially all the
results we shall describe in the present paper for conformal vertex
algebras also hold more generally for M\"obius vertex
algebras. However, for brevity, in the present paper we shall state
our results only for conformal vertex algebras. }
\end{rema}

As expected, a {\em module} for a conformal vertex algebra $V$ is
defined to be a module for $V$ viewed as a vertex algebra such that
the conformal element acts in the same way as in the definition of
vertex operator algebra.

In this paper we will study a larger family of ``modules'' which
satisfy all conditions for being a module except that they are only
direct sums of {\it generalized} eigenspaces, rather than eigenspaces,
for the operator $L(0)$.  We call these {\em generalized modules}.

We also have the obvious notions of {\em module homomorphism}, {\em
submodule}, {\em quotient module}, and so on.

For a generalized module $W$, we will denote by $W_{[n]}$ the
generalized eigenspace for $L(0)$ with respect to eigenvalue $n$; we
call $n$ the corresponding {\em (generalized) weight}.  Thus
\[
W=\coprod_{n\in\mathbb{C}} W_{[n]};
\]
and for $n\in\mathbb{C}$,
\[
W_{[n]}=\{w\in W\;|\; (L(0)-n)^m w=0 \;\mbox{ for }\;m\in
\mathbb{N}\;\mbox{ sufficiently large}\}.
\]

As in \cite{FHL}, from the $L(-1)$-derivative property and the Jacobi
identity for generalized modules for a conformal vertex algebra $V$, we
have
\[
[L(0),v_n]=(L(0)v)_n+(-n-1)v_n,
\]
as operators acting on a generalized $V$-module, for any $v\in V$ and
$n\in\mathbb{C}$. From this we obtain:
\begin{propo}\label{set:L(0)s}
Let $W$ be a generalized module for a conformal vertex algebra $V$.
Let $L(0)_s\in \mbox{End}\ W$ be the ``semisimple part'' of the
operator $L(0)$:
\[
L(0)_s w=n w \;\; \mbox{ for } \; w\in W_{[n]},\; n \in \mathbb{C}.
\]
Then the ``locally nilpotent part'' $L(0)-L(0)_s$ of $L(0)$ commutes
with the action of $V$ on $W$.  In other words, $L(0)-L(0)_s$ is a
$V$-homomorphism from $W$ to itself.
\end{propo}

Let $W$ be a generalized module for a conformal vertex algebra. We
will use the notation $\overline{W}$ for the completion of $W$ with
respect to the $\mathbb{C}$-grading, that is,
\[
\overline{W}=\prod_{n\in\mathbb{C}} W_{[n]}.
\]
We write
\begin{equation}\label{pin}
\pi_n: \overline{W}\to W_{[n]}
\end{equation}
for the natural projection map from $\overline{W}$ to $W_{[n]}$, for
any $n\in\mathbb{C}$.  We will also use the same notation for its
restriction to $W$, and for its natural extensions to spaces of formal
series with coefficients in $W$.

Clearly there is a canonical pairing between $\overline{W}$ and the
subspace $\coprod _{n\in \mathbb{C}} (W_{[n]})^*$ of $W^*$ (viewing
$(W_{[n]})^*$ as embedded in $W^*$ in the natural way).  Here and
throughout this paper we use the notation $M^*$ for the dual of a
vector space $M$. We will also use $\langle\cdot,\cdot\rangle$ for the
natural pairings between various vector spaces.

The following notion of ``strong gradedness'' for a conformal vertex
algebra will be needed:
\begin{defi}\label{def:dgv}
{\rm Let $A$ be an abelian group.  A conformal vertex algebra
\[
V=\coprod_{n\in \mathbb{Z}} V_{(n)}
\]
is said to be {\em strongly $A$-graded}, or just {\em strongly graded}
if $A$ is implicitly understood, if there is a second gradation on
$V$, by $A$,
\[
V=\coprod _{\alpha \in A} V^{(\alpha)},
\]
such that the following conditions are satisfied: the two gradations
are compatible, that is,
\[
V^{(\alpha)}=\coprod_{n\in \mathbb{Z}} V^{(\alpha)}_{(n)} \;\;
(\mbox{where}\;V^{(\alpha)}_{(n)}=V_{(n)}\bigcap V^{(\alpha)})\;
\mbox{ for any }\;\alpha \in A;
\]
for any $\alpha,\beta\in A$ and $n\in \mathbb{Z}$,
\begin{eqnarray*}
&V^{(\alpha)}_{(n)}=0\;\;\mbox{ for }\;n\;\mbox{ sufficiently
negative};&\\
&\dim V^{(\alpha)}_{(n)} <\infty;&\\
&{\bf 1}\in V^{(0)}_{(0)};&\\
&v_l V^{(\beta)} \subset V^{(\alpha+\beta)}\;\; \mbox{ for any }\;v\in
V^{(\alpha)},\;l\in \mathbb{Z};&
\end{eqnarray*}
and
\[
\omega\in V^{(0)}_{(2)}.
\]
}
\end{defi}

Given a strongly $A$-graded conformal vertex algebra $V$ and an
abelian group $\tilde A$ containing $A$ as a subgroup, a (generalized)
$V$-module $W$ is said to be {\em strongly $\tilde A$-graded}, or just
{\em strongly graded} if $\tilde A$ is implicitly understood, if the
obvious analogues of the conditions for $W$ in Definition
\ref{def:dgv} hold, except that the grading-truncation condition
becomes: for any $\beta\in \tilde A$ and $n\in \mathbb{C}$,
\begin{eqnarray}\label{set:dmltc}
&W^{(\beta)}_{[n+k]}=0 \;\;\mbox{ for }\;k\in \mathbb{Z}\;\mbox{
sufficiently negative.}&
\end{eqnarray}

\begin{rema}\label{trivialgroup}{\rm
Clearly, vertex operator algebras and their (ordinary) modules are
exactly the conformal vertex algebras and their (ordinary) modules
that are strongly graded with respect to the trivial group (with
(\ref{set:dmltc}) used as the grading-truncation condition).
Important examples of other strongly graded conformal vertex algebras
and modules come from vertex algebras and modules associated with an
even lattice that is not necessarily positive definite. We refer the
reader to \cite{Z} for the application of the results in the present
work to these examples. }
\end{rema}

Suppose that $L(1)$ acts locally nilpotently on a conformal vertex
algebra $V$ (which occurs, in particular, if $V$ is strongly
graded). Then generalizing formula (3.20) in \cite{tensor1} we define
the {\it opposite vertex operator} on a generalized $V$-module
$(W,Y_W)$ associated to $v\in V$ by
\begin{equation}\label{yo}
Y^o_W(v,x)=\sum_{n \in \mathbb{Z}} v^o_n
x^{-n-1}=Y_W(e^{x L(1)}(-x^{-2})^{L(0)}v,x^{-1}),
\end{equation}
so that
\[
v^o_n=(-1)^k\sum_{m\in \mathbb{N}}\frac{1}{m!}(L(1)^m v)_{-n-m-2+2k}
\]
for $v\in V_{(k)}$ and $n,k\in \mathbb{Z}$. (Note that the $L(1)$-local
nilpotence is needed for well-definedness here.) We have
\[
v^o_n W_{[m]}\subset W_{[m+n+1-k]}\;\;\mbox{ for any }\;m\;\in\mathbb{C},
\]
and as in the case of \cite{tensor1}, the opposite Jacobi identity
for $Y^o_{W}$, as in formula (3.23) of \cite{tensor1}, holds.

As in Section 5.2 of \cite{FHL}, we can define a $V$-action on $W^*$
as follows:
\begin{equation}\label{y'}
\langle Y'(v,x)w',w\rangle = \langle w', Y^o_W(v,x)w\rangle
\end{equation}
for $v\in V$, $w'\in W^*$ and $w\in W$; the correspondence $v\mapsto
Y'(v,x)$ is a linear map from $V$ to $({\rm End}\,W^*)[[x,x^{-1}]]$.
Define the operators $L'(n)\;(n\in \mathbb{Z})$ by
\[
Y'(\omega,x)=\sum_{n\in \mathbb{Z}}L'(n)x^{-n-2}.
\]

Now we can introduce the notion of ``contragredient module'' of a
strongly graded generalized module. First we give:
\begin{defi}{\rm
Let $W=\coprod_{\beta \in \tilde A,\;n\in\mathbb{C}}W^{(\beta)}_{[n]}$
be a strongly $\tilde A$-graded generalized module for a strongly
$A$-graded conformal vertex algebra.  We define $W'$ to be the
$(\tilde A\times\mathbb{C})$-graded vector space
\begin{equation}
W'=\coprod_{\beta \in \tilde A,\; n\in\mathbb{C}} (W')^{(\beta)}_{[n]},
\;\;\mbox{ where }\; (W')^{(\beta)}_{[n]}=(W^{(-\beta)}_{[n]})^*;
\end{equation}
also define
\[
(W')^{(\beta)}=\coprod_{n\in\mathbb{C}}(W^{(-\beta)}_{[n]})^*
\subset(W^{(-\beta)})^*
\]
and
\[
(W')_{[n]}=\coprod_{\beta\in\tilde A}(W^{(-\beta)}_{[n]})^*
\subset(W_{[n]})^*,
\]
the homogeneous subspaces of $W'$ with respect to the $\tilde A$- and
$\mathbb{C}$-grading, respectively.  There is of course a canonical
pairing between $W'$ and $\overline{W}\subset \prod_{\beta \in \tilde
A,\; n\in\mathbb{C}}W^{(\beta)}_{[n]}$.  }
\end{defi}

Let $W$ and $W'$ be as in this definition. It is straightforward to
show that the lower truncation condition for $Y'$ on $W'$ holds; thus
the Jacobi identity can be formulated on $W'$. Furthermore, using the
same arguments as in the proofs of Theorems 5.2.1 and 5.3.1 in
\cite{FHL} we have:
\begin{theo}\label{set:W'}
Let $\tilde A$ be an abelian group containing $A$ as a subgroup and
$V$ a strongly $A$-graded conformal vertex algebra. Let $(W,Y)$ be a
strongly $\tilde A$-graded $V$-module (respectively, generalized
$V$-module). Then the pair $(W',Y')$ carries a strongly $\tilde
A$-graded $V$-module (respectively, generalized $V$-module) structure,
and $(W'',Y'')=(W,Y)$.
\end{theo}

The pair $(W',Y')$ in Theorem \ref{set:W'} will be called the {\em
contragredient module} of $(W,Y)$.

{\it Throughout this paper, we will fix a conformal vertex algebra
$V$.}  Given abelian groups $A$ and $\tilde A \supset A$ as above,
when $V$ is strongly $A$-graded we will write ${\cal C}_1$ for the
category of all strongly $\tilde A$-graded generalized $V$-modules.
We have a contravariant functor $(\cdot)': (W,Y)\mapsto (W',Y')$ from
${\cal C}_1$ to itself, which we call the {\em contragredient
functor}.  Our main objects of study will be certain full
subcategories ${\cal C}$ of ${\cal C}_1$ that are closed under the
contragredient functor.

\setcounter{equation}{0}
\setcounter{rema}{0}

\section{Logarithmic intertwining operators}
In this section we first introduce the ``logarithm of a formal
variable'' and study some of its properties. We then study the notion
of ``logarithmic intertwining operator'' introduced in \cite{Mi} and
give some results which are essential for our generalization of the
tensor product theory.

We will use the following usual notation for formal series with
arbitrary complex powers, as in \cite{FLM}:
\[
{\cal W}\{x\}=\bigg\{\sum_{n\in\mathbb{C}}w(n)x^n\,|\,w(n)\in {\cal W},
\;n\in\mathbb{C}\bigg\}
\]
for any vector space ${\cal W}$ and formal variable $x$.

{}From now on we will sometimes need and use new independent
(commuting) formal variables denoted by $\log x$, $\log y$, $\log x_1$,
$\log x_2, \dots,$ etc. We will work with formal series in such formal
variables together with the ``usual'' formal variables $x$, $y$,
$x_1$, $x_2, \dots,$ etc., with coefficients in certain vector spaces,
and the powers of the monomials in {\it all} the variables can be
arbitrary complex numbers.  (Later we will restrict our attention to
only nonnegative integral powers of the ``$\log$'' variables.)

Given a formal variable $x$, we write $\frac d{dx}$ for the linear map
on ${\cal W}\{x,\log x\}$, for any vector space ${\cal W}$ not
involving $x$, defined (and indeed well defined) by the (expected)
formula
\begin{eqnarray*}
\lefteqn{\frac d{dx}\bigg(\sum_{m,n\in \mathbb{C}}w_{n,m}x^n(\log
x)^m\bigg)}\\
&&\ =\sum_{m,n\in\mathbb{C}}((n+1)w_{n+1,m}+ (m+1)w_{n+1,m+1})x^{n}(\log
x)^{m}\\
&&\bigg(=\sum_{m,n\in\mathbb{C}}n w_{n,m}x^{n-1}(\log x)^m+ \sum_{m,n\in
\mathbb{C}}m w_{n,m}x^{n-1}(\log x)^{m-1}\bigg)
\end{eqnarray*}
where $w_{n,m}\in {\cal W}$ for all $m, n\in \mathbb{C}$. The same
notation will also be used for the restriction of $\frac d{dx}$ to any
subspace of ${\cal W}\{x,\log x\}$ that is closed under $\frac d{dx}$,
e.g., ${\cal W}\{x\}[[\log x]]$ or $\mathbb{C}[x,x^{-1},\log
x]$. Clearly, $\frac d{dx}$ acting on ${\cal W}\{x\}$ coincides with
the usual formal derivative.

\begin{rema}{\rm
Let $f$, $g$ and $f_i$, $i$ in some index set $I$, all be of the form
\begin{equation}\label{log:f}
\sum_{m,n\in \mathbb{C}}w_{n,m}x^n(\log x)^m\in {\cal W}\{x,\log
x\},\;\; w_{n,m}\in {\cal W}.
\end{equation}
One checks the following straightforwardly: Suppose that the sum of
$f_i$, $i\in I$, exists (in the obvious sense). Then the sum of
$\frac d{dx}f_i$, $i\in I$, also exists and is equal to $\frac
d{dx}\sum_{i\in I}f_i$.  More generally, for any $T=p(x)\frac{d}{dx}$,
$p(x)\in \mathbb{C}[x,x^{-1}]$, the sum of $Tf_i$, $i\in I$, exists and
is equal to $T\sum_{i\in I}f_i$. Thus the sum of $e^{yT}f_i$, $i\in
I$, exists and is equal to $e^{yT}\sum_{i\in I}f_i$ ($e^{yT}$ being
the formal exponential series, as usual).  Suppose that ${\cal W}$ is
an (associative) algebra or that the coefficients of either $f$ or $g$
are complex numbers. If the product of $f$ and $g$ exists, then the
product of $\frac d{dx}f$ and $g$ and the product of $f$ and $\frac
d{dx}g$ both exist, and $\frac d{dx}(fg)=(\frac d{dx}f)g+f(\frac
d{dx}g)$. Furthermore, for any $T$ as before, the product of $Tf$ and
$g$ and the product of $f$ and $Tg$ both exist, and
$T(fg)=(Tf)g+f(Tg)$. In addition, the product of $e^{yT}f$ and
$e^{yT}g$ exists and is equal to $e^{yT}(fg)$, just as in formulas
(8.2.6)--(8.2.10) of \cite{FLM}.  The point here, of course, is just
the formal derivation property of $\frac d{dx}$, except that sums and
products of expressions do not exist in general.  }
\end{rema}

\begin{rema}{\rm
Note that the ``equality'' $x=e^{\log x}$ does not hold, since the
left-hand side is a formal variable, while the right-hand side is a
formal series in another formal variable. In fact, this formula should
not be assumed, since, for example, the formal delta function
$\delta(x)=\sum_{n\in \mathbb{Z}}x^n$ would not exist in the sense of
formal calculus, if $x$ were allowed to be replaced by the formal
series $e^{\log x}$. By contrast, note that the equality
\begin{equation}\label{log:logex}
\log e^x=x
\end{equation}
does indeed hold. This is because the formal series $e^x$ is of the
form $1+X$ where $X$ involves only positive integral powers of $x$ and
in (\ref{log:logex}), ``$\log$'' refers to the usual formal
logarithmic series
\[
\log(1+X)=\sum_{i\geq 1} \frac{(-1)^{i-1}}i X^i,
\]
{\em not} to the ``$\log$'' of a formal variable. We will use the
symbol ``$\log$'' in both ways, and the meaning will be clear in
context.  }
\end{rema}

We will typically use notations of the form $f(x)$, instead of
$f(x,\log x)$, to denote elements of ${\cal W}\{x,\log x\}$ for some
vector space ${\cal W}$ as above. For this reason, we need to
interpret the meaning of symbols such as $f(x+y)$, or more generally,
symbols obtained by replacing $x$ in $f(x)$ by something other than
just a single formal variable (since $\log x$ is a formal variable and
not the image of some operator acting on $x$).  Specifically, we use
the following notational conventions; the existence of the expressions
will be justified in Remark \ref{log:exist}:
\begin{nota}{\rm
For formal variables $x$, $y$, and $f(x)$ of the form (\ref{log:f}),
define
\begin{eqnarray}\label{log:not1}
f(x+y)&=&\sum_{m,n\in \mathbb{C}}w_{n,m}(x+y)^n\bigg(\log x+
\log\bigg(1+\frac yx\bigg)\bigg)^m\nno\\
&=&\sum_{m,n\in \mathbb{C}}w_{n,m}(x+y)^n\bigg(\log x+
\sum_{i\geq 1} \frac{(-1)^{i-1}}i \bigg(\frac yx\bigg)^i\bigg)^m;
\end{eqnarray}
in the right-hand side, $(\log x+\sum_{i\geq 1}\frac{(-1)^{i-1}}i
(\frac yx)^i)^m$, according to the binomial expansion convention, is
to be expanded in nonnegative integral powers of the second summand
$\sum_{i\geq 1} \frac{(-1)^{i-1}}i (\frac yx)^i$, so the right-hand
side of (\ref{log:not1}) is equal to
\begin{equation}\label{log:1-tmp}
\sum_{m,n\in \mathbb{C}}w_{n,m}(x+y)^n\sum_{j\in \mathbb{N}}{m\choose j}
(\log x)^{m-j}\bigg(\sum_{i\geq 1} \frac{(-1)^{i-1}}i \bigg(\frac
yx\bigg)^i\bigg)^j
\end{equation}
when expanded one step further. Also define
\begin{equation}\label{log:not2}
f(xe^y)=\sum_{m,n\in \mathbb{C}}w_{n,m}x^ne^{ny}(\log x+y)^m,
\end{equation}
\begin{equation}\label{log:not3}
f(xy)=\sum_{m,n\in \mathbb{C}}w_{n,m}x^ny^n(\log x+\log y)^m.
\end{equation}
}
\end{nota}

\begin{rema}\label{log:exist}{\rm
The existence of the right-hand side of (\ref{log:not1}), or
(\ref{log:1-tmp}), can be seen by writing $(x+y)^n$ as $x^n(1+\frac
yx)^n$ and observing that

\[
\bigg(\sum_{i\geq 1} \frac{(-1)^{i-1}}i \bigg(\frac yx\bigg)^i\bigg)^j
\in \bigg(\frac yx\bigg)^j\mathbb{C}\bigg[\bigg[\frac yx\bigg]\bigg].
\]
The existence of the right-hand sides of (\ref{log:not2}) and of
(\ref{log:not3}) is clear. Furthermore, both $f(x+y)$ and $f(xe^y)$
lie in ${\cal W}\{x,\log x\}[[y]]$, while $f(xy)$ lies in ${\cal
W}\{xy,\log x\}[[\log y]]$.  }
\end{rema}

We still have the following identities (cf. Section 8.3 of
\cite{FLM}) in the logarithmic settings:
\begin{propo}\label{log:ck12}
For $f(x)$ as in (\ref{log:f}), we have
\begin{equation}\label{log:ck1}
e^{y\frac d{dx}}f(x)=f(x+y)
\end{equation}
and
\begin{equation}\label{log:ck2}
e^{yx\frac d{dx}}f(x)=f(xe^y).
\end{equation}
\end{propo}

While (\ref{log:ck2}) is straightforward to verify, (\ref{log:ck1}) is
surprisingly subtle (when logarithmic terms are involved). For
instance, it turns out that (\ref{log:ck1}) essentially amounts to a
generating function form of the following combinatorial identity:
\begin{eqnarray}\label{log:comb}
\frac{j!}{k!}\sum_{0<t_1<t_2<\cdots<t_{k-j}<k} t_1t_2\cdots t_{k-j}=
\sum_{\tiny\begin{array}{c}i_1+\cdots+i_j=k\\
1\leq i_1,\dots,i_j\leq k\end{array}} \frac{1}{i_1i_2\cdots i_j}
\end{eqnarray}
for all $k\in \mathbb{N}$ and $j=0,\dots,k$.  A wide variety of
identities of this type will clearly arise as a result of the deeper
part of our theory.  Incidentally, at the time this work was being
written, one of us (L.Z.) happened to pick up the then-current issue
of the American Mathematical Monthly and happened to notice the
following problem {}from the Problems and Solutions section, proposed
by D. Lubell \cite{Lu}:
\begin{quote}
Let $N$ and $j$ be positive integers, and let $S=\{(w_1,\dots, w_j)\in
\mathbb{Z}_+^j\,|\,0<w_1+\cdots+w_j\leq N\}$ and
$T=\{(w_1,\cdots,w_j)\in \mathbb{Z}_+^j\,|\linebreak w_1,\dots,w_j
\mbox{ are distinct and bounded by }N\}.$ Show that
$$\sum_S\frac 1{w_1\cdots w_j}=\sum_T\frac 1{w_1\cdots w_j}.$$
\end{quote}
This follows immediately from (\ref{log:comb}), which is in fact a
refinement, since the left-hand side of (\ref{log:comb}) is equal to
\[
j!\sum_{1\leq w_1<w_2<\cdots<w_{j-1}\leq k-1} \frac 1{w_1w_2\cdots
w_{j-1}k}=\sum_{T_k} \frac 1{w_1w_2\cdots w_j}
\]
where $T_k=\{(w_1,\dots,w_j)\in\{1,2,\dots,k\}^j\,|\,w_i\mbox{
distinct, with maximum exactly }k\}$, the right-hand side is
\[
\sum_{S_k} \frac 1{w_1w_2\cdots w_j}
\]
where $S_k=\{(w_1,\dots,w_j)\in\{1,2,\dots,k\}^j\,|\,
w_1+\cdots+w_j=k\}$, and one has $S=\coprod_{k=1}^N S_k$ and
$T=\coprod_{k=1}^N T_k$.

When we define the notion of logarithmic intertwining operator below,
we will impose a condition requiring certain formal series to lie in
spaces of the type ${\cal W}[\log x]\{x\}$ (so that for each power of
$x$, possibly complex, we have a {\it polynomial} in $\log x$), partly
because such results as the following (which is expected) will indeed
hold in our formal setup when the powers of the formal variables are
restricted in this way.
\begin{lemma}
Let $a\in \mathbb{C}$ and $m\in \mathbb{Z}_+$. If $f(x)\in {\cal W}[\log
x]\{x\}$ (${\cal W}$ any vector space not involving $x$ or $\log x$)
satisfies the formal differential equation
\begin{equation}\label{de:(xdx-a)^m}
\bigg(x\frac{d}{dx}-a\bigg)^mf(x)=0,
\end{equation}
then $f(x)\in {\cal W}x^a\oplus{\cal W}x^a\log x \cdots\oplus{\cal
W}x^a(\log x)^{m-1}$; and furthermore, if $m$ is the smallest integer
so that (\ref{de:(xdx-a)^m}) is satisfied, then the coefficient of
$x^a(\log x)^{m-1}$ in $f(x)$ is nonzero.
\end{lemma}

This can be proved by induction on $m$.

Note that there are of course solutions of the equation
(\ref{de:(xdx-a)^m}) outside the space ${\cal W}[\log x]\{x\}$, for
example, $f(x)=wx^be^{(a-b)\log x}\in x^b{\cal W}[[\log x]]$ for any
complex number $b \neq a$ and any nonzero vector $w\in {\cal W}$.

Following \cite{Mi}, with a slight generalization, we now introduce
the notion of logarithmic intertwining operator. We will later see
that the axioms in the following definition correspond exactly to
those in the notion of certain ``intertwining maps'' (cf. Definition
\ref{im:imdef}).

\begin{defi}\label{def:logio}{\rm
Let $W_1$, $W_2$, $W_3$ be generalized $V$-modules.
A {\em logarithmic intertwining operator of type
${W_3\choose W_1\,W_2}$} is a linear map
\begin{equation}\label{log:map0}
{\cal Y}(\cdot, x)\cdot: W_1\otimes W_2\to W_3[\log x]\{x\},
\end{equation}
or equivalently,
\begin{equation}\label{log:map}
w_{(1)}\otimes w_{(2)}\mapsto{\cal Y}(w_{(1)},x)w_{(2)}=\sum_{n\in
\mathbb{C}}\sum_{k\in \mathbb{N}}{w_{(1)}}_{n;\,k}^{\cal
Y}w_{(2)}x^{-n-1}(\log x)^k\in W_3[\log x]\{x\}
\end{equation}
for all $w_{(1)}\in W_1$ and $w_{(2)}\in W_2$, such that the following
conditions are satisfied: the {\em lower truncation condition}: for
any $w_{(1)}\in W_1$, $w_{(2)}\in W_2$ and $n\in \mathbb{C}$,
\begin{equation}\label{log:ltc}
{w_{(1)}}_{n+m;\,k}^{\cal Y}w_{(2)}=0\;\;\mbox{ for }\;m\in \mathbb{N}
\;\mbox{ sufficiently large,\, independently of}\;k;
\end{equation}
the {\em Jacobi identity}:
\begin{eqnarray}\label{log:jacobi}
\lefteqn{\dps x^{-1}_0\delta \bigg( {x_1-x_2\over x_0}\bigg)
Y(v,x_1){\cal Y}(w_{(1)},x_2)w_{(2)}}\nno\\
&&\hspace{2em}- x^{-1}_0\delta \bigg( {x_2-x_1\over -x_0}\bigg)
{\cal Y}(w_{(1)},x_2)Y(v,x_1)w_{(2)}\nno \\
&&{\dps = x^{-1}_2\delta \bigg( {x_1-x_0\over x_2}\bigg) {\cal
Y}(Y(v,x_0)w_{(1)},x_2) w_{(2)}}
\end{eqnarray}
for $v\in V$, $w_{(1)}\in W_1$ and $w_{(2)}\in W_2$ (note that the
first term on the left-hand side is meaningful because of
(\ref{log:ltc})); and the {\em $L(-1)$-derivative property:} for any
$w_{(1)}\in W_1$,
\begin{equation}\label{log:L(-1)dev}
{\cal Y}(L(-1)w_{(1)},x)=\frac d{dx}{\cal Y}(w_{(1)},x).
\end{equation}
Suppose in addition that
$V$ and $W_1$, $W_2$ and $W_3$ are strongly graded.
A logarithmic intertwining operator
${\cal Y}$ is a {\em grading-compatible
logarithmic intertwining operator} if for $\beta, \gamma \in \tilde A$
and $w_{(1)} \in W_{1}^{(\beta)}$,
$w_{(2)} \in W_{2}^{(\gamma)}$, $n \in \mathbb{C}$ and $k \in \mathbb{N}$, we
have
\begin{equation}\label{gradingcompatcondn}
{w_{(1)}}_{n;\,k}^{\cal Y}w_{(2)} \in W_{3}^{(\beta + \gamma)}.
\end{equation}}
\end{defi}

By setting $v=\omega$ in (\ref{log:jacobi}) and then taking
$\res_{x_0}\res_{x_1}x_1^{j+1}$ we have
\begin{equation}\label{log:L(j)b}
{}[L(j), {\cal Y}(w_{(1)},x)]=\sum_{i=0}^{j+1}{j+1\choose i}x^i{\cal
Y}(L(j-i)w_{(1)},x)
\end{equation}
for any $w_{(1)}\in W_1$ and $j = -1,0,1$.

Every ordinary intertwining operator (as in, for example,
\cite{tensor1}) is of course a logarithmic intertwining operator which
involves no formal variable $\log x$.  (The present lower truncation
condition is more relaxed than the one in \cite{tensor1}.)  The
grading-compatible logarithmic intertwining operators of a fixed type
${W_3\choose W_1\, W_2}$ form a vector space, which we denote by
${\cal V}^{W_3}_{W_1\,W_2}$.  We call the dimension of ${\cal
V}^{W_3}_{W_1\,W_2}$ the {\it fusion rule} for $W_1$, $W_2$ and $W_3$
and denote it by $N^{W_3}_{W_1\,W_2}$ (see Remark \ref{log:ordi}
below).

By direct analysis we can prove:
\begin{lemma}\label{log:lemma}
Let $W_1$, $W_2$, $W_3$ be generalized $V$-modules. Let
\begin{eqnarray*}
{\cal Y}(\cdot,x)\cdot: W_1\otimes W_2 & \to & W_3\{x,\log x\}\\
w_{(1)}\otimes w_{(2)}&\mapsto&{\cal Y}(w_{(1)},x)w_{(2)}=
\sum_{n,k\in \mathbb{C}}{w_{(1)}}_{n;\,k}^{\cal Y}w_{(2)}x^{-n-1}(\log
x)^k
\end{eqnarray*}
be a linear map that satisfies the L(-1)-derivative property
(\ref{log:L(-1)dev}) and the $L(0)$-bracket relation, that is,
(\ref{log:L(j)b}) with $j=0$. Then for any $a,b,c\in\mathbb{C}$,
$t\in\mathbb{N}$, $w_{(1)}\in W_1$ and $w_{(2)}\in W_2$,
\begin{eqnarray*}
\lefteqn{(L(0)-c)^t{\cal Y}(w_{(1)},x)w_{(2)}=\sum_{i,j,l\in
\mathbb{N},\;i+j+l=t}\frac{t!}{i!j!l!}\cdot}\\
&&\cdot\left(x\frac d{dx}-c+a+b\right)^l{\cal Y}((L(0)-a)^iw_{(1)},x)
(L(0)-b)^jw_{(2)}.
\end{eqnarray*}
Also, for any $a,b,n,k\in\mathbb{C}$, $t\in\mathbb{N}$, $w_{(1)}\in W_1$
and $w_{(2)}\in W_2$, we have
\begin{eqnarray*}
\lefteqn{(L(0)-a-b+n+1)^t({w_{(1)}}_{n;\,k}^{\cal Y}w_{(2)})}\\
&&=t!\sum_{i,j,l\geq 0,\;i+j+l=t}{k+l\choose l}
\bigg(\frac{(L(0)-a)^i}{i!}w_{(1)}\bigg)^{\cal
Y}_{n;\,k+l}\bigg(\frac{(L(0)-b)^j}{j!}w_{(2)}\bigg);
\end{eqnarray*}
in generating function form, this gives
\begin{eqnarray}\label{lmaf-4}
\lefteqn{e^{y(L(0)-a-b+n+1)}({w_{(1)}}_{n;\,k}^{\cal
Y}w_{(2)})}\nno\\
&&=\sum_{l\in\mathbb{N}}{k+l\choose
l}(e^{y(L(0)-a)}w_{(1)})_{n;\,k+l}^{\cal Y}(e^{y(L(0)-b)}w_{(2)})y^l.
\end{eqnarray}
\end{lemma}

Using this lemma, we have the following result, which summarizes some
important properties of logarithmic intertwining operators:
\begin{propo}\label{log:logwt}
Let $W_1$, $W_2$, $W_3$ be generalized $V$-modules, and let ${\cal
Y}(\cdot, x)\cdot$ be a logarithmic intertwining operator of type
${W_3\choose W_1\,W_2}$.  Let $w_{(1)}$ and $w_{(2)}$ be homogeneous
elements of $W_1$ and $W_2$ of (generalized) weights $n_1$ and $n_2
\in \mathbb{C}$, respectively, and let $k_1$ and $k_2$ be positive
integers such that $(L(0)-n_1)^{k_1}w_{(1)}=0$ and
$(L(0)-n_2)^{k_2}w_{(2)}=0$.  Then we have:

(a) (\cite{Mi}) For any $w'_{(3)}\in W_3^*$, $n_3\in \mathbb{C}$ and
$k_3\in \mathbb{Z}_+$ such that $(L'(0)-n_3)^{k_3}w'_{(3)}=0$,
\begin{eqnarray*}
\lefteqn{\langle w'_{(3)}, {\cal Y}(w_{(1)}, x)w_{(2)}\rangle}\\
&&\in \mathbb{C}x^{n_3-n_1-n_2}\oplus\mathbb{C}x^{n_3-n_1-n_2}\log
x\oplus\cdots\oplus \mathbb{C}x^{n_3-n_1-n_2}(\log
x)^{k_1+k_2+k_3-3}.
\end{eqnarray*}

(b) For any $n\in \mathbb{C}$ and $k\in \mathbb{N}$, ${w_{(1)}}^{\cal
Y}_{n;\,k}w_{(2)}\in W_3$ is homogeneous of (generalized) weight
$n_1+n_2-n-1$.

(c) Fix $n\in \mathbb{C}$ and $k\in \mathbb{N}$. For each $i,j\in
\mathbb{N}$, let $m_{ij}$ be a nonnegative integer such that
\[
(L(0)-n_1-n_2+n+1)^{m_{ij}}(((L(0)-n_1)^iw_{(1)})^{\cal
Y}_{n;\,k}(L(0)-n_2)^jw_{(2)})=0.
\]
Then for all $t\geq \max\{m_{ij}\,|\,0\leq i<k_1,\; 0\leq
j<k_2\}+k_1+k_2-2$,
\[
{w_{(1)}}^{\cal Y}_{n;\,k+t}w_{(2)}=0.
\]
\end{propo}

Proposition \ref{log:logwt}(b) immediately gives:
\begin{corol}
Let $W_1$, $W_2$ and $W_3$
be generalized $V$-modules whose weights are all congruent modulo
$\mathbb{Z}$ to complex numbers $h_1$, $h_2$ and $h_3$, respectively
(for example, when $W_1$, $W_2$ and $W_3$ are all indecomposable).
Let ${\cal Y}(\cdot, x)\cdot$ be a logarithmic intertwining operator
of type ${W_3\choose W_1\,W_2}$.  Then all the powers of $x$ in ${\cal
Y}(\cdot, x)\cdot$ are congruent modulo $\mathbb{Z}$ to $h_3-h_1-h_2$.
\end{corol}

\begin{rema}\label{log:ordi}{\rm
Let $W_1$, $W_2$ and $W_3$ be (ordinary) $V$-modules.  Then it follows
{}from Proposition \ref{log:logwt}(a), or alternatively, {}from
Proposition \ref{log:logwt}(b) and (c), that any logarithmic
intertwining operator of type ${W_3\choose W_1\,W_2}$ is just an
ordinary intertwining operator of this type, i.e., it does not involve
$\log x$.  As a result, when $V$ is a vertex operator algebra, for
(ordinary) $V$-modules the notion of fusion rule defined in this paper
coincides with the notion of fusion rule defined in, for example,
\cite{tensor1} (except for the minor issue of the truncation
condition.)}
\end{rema}

Our definition of logarithmic intertwining operator is identical to
that in \cite{Mi} (in case $V$ is a vertex operator algebra) except
that in \cite{Mi}, a logarithmic intertwining operator ${\cal Y}$ of
type ${W_3\choose W_1\,W_2}$ is required to be a linear map $W_1\to
\hom(W_2,W_3)\{x\}[\log x]$, instead of as in (\ref{log:map0}), and
the lower truncation condition (\ref{log:ltc}) is replaced by:
\[
{w_{(1)}}^{\cal Y}_{n;k}w_{(2)}=0\;\;\mbox{for}\;n\;\mbox{whose real
part is sufficiently large}
\]
for any $w_{(1)}\in W_1$, $w_{(2)}\in W_2$ and $k\in \mathbb{N}$.

Given a logarithmic intertwining operator ${\cal Y}$ as in
(\ref{log:map}), set
\[
{\cal Y}^{(k)}(w_{(1)},x)w_{(2)}=\sum_{n\in
\mathbb{C}}{w_{(1)}}_{n;\,k}^{\cal Y}w_{(2)}x^{-n-1}
\]
for $k\in \mathbb{N}$, $w_{(1)}\in W_1$ and $w_{(2)}\in W_2$, so that
\[
{\cal Y}(w_{(1)},x)w_{(2)}=\sum_{k\in \mathbb{N}}{\cal Y}^{(k)}
(w_{(1)},x)w_{(2)}(\log x)^k.
\]
Then each ${\cal Y}^{(k)}$ satisfies the Jacobi identity, but not in
general the $L(-1)$-derivative property. However, the following result
shows that suitable formal linear combinations of certain
modifications of ${\cal Y}^{(k)}$ (depending on $t\in\mathbb{N}$; see
below) form a sequence of logarithmic intertwining operators:
\begin{propo}
Let $W_1$, $W_2$, $W_3$ be generalized $V$-modules,
and let ${\cal Y}(\cdot, x)\cdot$ be a logarithmic
intertwining operator of type ${W_3\choose W_1\,W_2}$.  For $\mu\in
\mathbb{C}/\mathbb{Z}$ and $t \in \mathbb{N}$, define ${\cal
X}^{\mu}_t:W_1\otimes W_2\to W_3[\log x]\{x\}$ by:
\[
{\cal X}^{\mu}_t: w_{(1)}\otimes w_{(2)}\mapsto \sum_{k\in \mathbb{N}}
{k+t\choose t}\sum_{\bar n=\mu}{w_{(1)}}_{n; \,k+t}^{\cal Y}
w_{(2)}x^{-n-1}(\log x)^k.
\]
Then each ${\cal X}^{\mu}_t$ is a logarithmic intertwining operator of
type ${W_3\choose W_1\,W_2}$. In particular, the operator ${\cal X}_t$
defined by
\begin{equation}\label{newio}
{\cal X}_t: w_{(1)}\otimes w_{(2)}\mapsto \sum_{k\in \mathbb{N}}
{k+t\choose t}\sum_{n\in\mathbb{C}}{w_{(1)}}_{n; \,k+t}^{\cal Y}
w_{(2)}x^{-n-1}(\log x)^k
\end{equation}
is a logarithmic intertwining operator of the same type.  In the
strongly graded case, if ${\cal Y}$ is grading-compatible, then so are
${\cal X}^{\mu}_t$ and ${\cal X}_t$.
\end{propo}

This can be proved essentially by induction on $t$.

\begin{rema}{\rm
Given any logarithmic intertwining operator ${\cal Y}(\cdot,x)\cdot$
as in (\ref{log:map0}) and any $i,j,k\in \mathbb{N}$, by Proposition
\ref{set:L(0)s} we see that
\[
(L(0)-L(0)_s)^k{\cal Y}((L(0)-L(0)_s)^i\cdot,x) (L(0)-L(0)_s)^j\cdot
\]
is again a logarithmic intertwining operator, and in the strongly
graded case, if ${\cal Y}$ is grading-compatible, so is this operator.
It turns out that the intertwining operators (\ref{newio}) are just
linear combinations of these.  }
\end{rema}

Now we define operators ``$x^{\pm L(0)}$'' for generalized modules, in
the natural way:
\begin{defi}{\rm
Let $W$ be a generalized $V$-module. We
define $x^{\pm L(0)}: W\to W\{x\}[\log x]\subset W[\log x]\{x\}$ as
follows: For any $w\in W_{[n]}$ ($n\in\mathbb{C}$), define
\[
x^{\pm L(0)}w=x^{\pm n}e^{\pm\log x(L(0)-n)}w
\]
(note that the local nilpotence of $L(0)-n$ on $W_{[n]}$ insures that
the formal exponential series terminates) and then extend linearly to
all $w\in W$. We also define operators $x^{\pm L'(0)}$ on $W^*$ by the
condition that for all $w'\in W^*$ and $w\in W$,
\[
\langle x^{\pm L'(0)}w',w\rangle=\langle w',x^{\pm L(0)}w\rangle\;
(\in \mathbb{C}\{x\}[\log x]),
\]
so that $x^{\pm L'(0)}: W^*\to W^*\{x\}[[\log x]]$.  }
\end{defi}

\begin{rema}{\rm
Note that for $w\in W_{[n]}$, by definition we have
\begin{equation}\label{log:x^L(0)}
x^{\pm L(0)}w=x^{\pm n}\sum_{i\in \mathbb{N}}\frac{(L(0)-n)^iw}{i!}(\pm
\log x)^i\in x^{\pm n}W_{[n]}[\log x].
\end{equation}
It is also easy to verify that for any $w\in W$,
\begin{eqnarray}
&\displaystyle
x^{L(0)}x^{-L(0)}w=w=x^{-L(0)}x^{L(0)}w,&\label{=w=}\\
&\displaystyle\frac{d}{dx}x^{\pm L(0)}w=\pm x^{-1}x^{\pm
L(0)}L(0)w.&
\end{eqnarray}
}
\end{rema}

We have the following generalizations to logarithmic intertwining
operators of three standard formulas for (ordinary) intertwining
operators (see \cite{FHL}, formulas (5.4.21), (5.4.22) and (5.4.23));
see also \cite{Mi} for parts (a) and (b):
\begin{propo}
Let ${\cal Y}$ be a logarithmic intertwining operator of type
${W_3\choose W_1\,W_2}$ and let $w\in W_1$. Then
\begin{description}
\item{(a)}
\[
e^{yL(-1)}{\cal Y}(w,x)e^{-yL(-1)}={\cal Y}(e^{yL(-1)}w,x)={\cal
Y}(w,x+y)
\]
(recall (\ref{log:not1}))
\item{(b)}
\[
y^{L(0)}{\cal Y}(w,x)y^{-L(0)}={\cal Y}(y^{L(0)}w,xy)
\]
(recall (\ref{log:not3}))
\item{(c)}
\[
e^{yL(1)}{\cal Y}(w,x)e^{-yL(1)}={\cal
Y}(e^{y(1-yx)L(1)}(1-yx)^{-2L(0)}w,x(1-yx)^{-1}).
\]
\end{description}
\end{propo}

The equality (a) essentially follows from the identity
\[
L(-1){\cal Y}(w,x)={\cal Y}(L(-1)w,x)+{\cal Y}(w,x)L(-1),\quad w\in W,
\]
and (\ref{log:ck1}) in Proposition \ref{log:ck12}, (b) essentially
follows from (\ref{log:not3}), (\ref{lmaf-4}) and (\ref{=w=}), and the
proof of (c) uses results in \cite{FHL}.

Let ${\cal Y}$ be a logarithmic intertwining operator of type
${W_3\choose W_1\, W_2}$. By analogy with the ordinary case in Section
7.1 of \cite{tensor2}, for any complex number $\zeta$, $w_{(1)}\in
W_1$ and $w_{(2)}\in W_2$,
\begin{equation}\label{log:subs}
{\cal Y}(w_{(1)},y)w_{(2)}\lbar_{y^n=e^{n\zeta}x^n,\;
(\log y)^k=(\zeta+\log x)^k,\;n\in\mathbb{C},\;k\in\mathbb{N}}
\end{equation}
is a well-defined element of $W_3[\log x]\{x\}$. We denote this
element by ${\cal Y} (w_{(1)},e^\zeta x)w_{(2)}$.  Note that it
depends on $\zeta$, not just on $e^\zeta$. Given any $r\in
\mathbb{Z}$, we define $\Omega_r({\cal Y}): W_2\otimes W_1\to W_3[\log
x]\{x\}$ by
\begin{equation}\label{Omega}
\Omega_{r}({\cal Y})(w_{(2)},x)w_{(1)} = e^{xL(-1)} {\cal
Y}(w_{(1)},e^{ (2r+1)\pi i}x)w_{(2)}
\end{equation}
for $w_{(1)}\in W_{1}$ and $w_{(2)}\in W_{2}$, the same formula as
formula (7.1) of \cite{tensor2}. We can prove that Proposition 7.1 in
\cite{tensor2} still holds when we replace the phrase ``an
intertwining operator ${\cal Y}$'' by ``a logarithmic intertwining
operator ${\cal Y}$.''  This result continues to hold for
grading-compatible logarithmic intertwining operators, in the strongly
graded case, and in particular, ${\cal Y}\mapsto \Omega_r({\cal Y})$
is a linear isomorphism from ${\cal V}_{W_1\,W_2}^{W_3}$ to ${\cal
V}_{W_2\,W_1}^{W_3}$ with inverse given by $\Omega_{-r-1}$.

In case $V$, $W_1$, $W_2$ and $W_3$ are strongly graded, we can also
carry over the concept of ``$r$-contragredient operator'' from
\cite{tensor2}, as follows: Given a grading-compatible logarithmic
intertwining operator ${\cal Y}$ of type ${W_3\choose W_1\,W_2}$ and
an integer $r$, we define the {\em $r$-contragredient operator of
${\cal Y}$} to be the linear map
\begin{eqnarray*}
W_1\otimes W'_3&\to&W'_2\{x\}[[\log x]]\\
w_{(1)}\otimes w'_{(3)}&\mapsto&A_r({\cal Y})(w_{(1)},x)w'_{(3)}
\end{eqnarray*}
given by
\begin{eqnarray*}
\lefteqn{\langle A_r({\cal Y})(w_{(1)},x)w'_{(3)},w_{(2)}
\rangle_{W_2}}\nno\\
&&=\langle w'_{(3)},{\cal Y}(e^{xL(1)}e^{(2r+1)\pi iL(0)}(x^{-L(0)})^2
w_{(1)},x^{-1})w_{(2)}\rangle_{W_3},
\end{eqnarray*}
for any $w_{(1)}\in W_1$, $w_{(2)}\in W_2$ and $w'_{(3)}\in W'_3$,
where we use the notation
\[
f(x^{-1})=\sum_{m\in\mathbb{N},\,n\in\mathbb{C}}w_{n,m} x^{-n}(-\log x)^m
\]
for any $f(x)=\sum_{m\in\mathbb{N},\,n\in \mathbb{C}}w_{n,m} x^n(\log
x)^m\in{\cal W}\{x\}[[\log x]]$, ${\cal W}$ any vector space (not
involving $x$).  Note that for the case $W_1=V$, $W_2=W_3=W$ and
${\cal Y}=Y_W$, the operator $A_r(Y_W)$ agrees with the contragredient
vertex operator $Y'_W$ of $Y_W$ (recall (\ref{yo}) and (\ref{y'})) for
any $r\in\mathbb{Z}$.  In general, we can prove that Proposition 7.3
in \cite{tensor2} still holds when we replace the phrase ``an
intertwining operator ${\cal Y}$'' by ``a grading-compatible
logarithmic intertwining operator ${\cal Y}$,'' and in particular
${\cal Y}\mapsto A_r({\cal Y})$ defines a linear isomorphism from
${\cal V}_{W_1\,W_2}^{W_3}$ to ${\cal V}_{W_1\,W'_3}^{W'_2}$, with
inverse given by $A_{-r-1}$.

Let $V$, $W_{1}$, $W_2$ and $W_3$ be strongly graded. Set
\[
N_{W_1W_2W_3}=N_{W_1\,W_2}^{W'_3}.
\]
Then from the above discussion, we see that, as in the ordinary case,
for any permutation $\sigma$ of $(1,2,3)$,
$N_{W_1W_2W_3}=N_{W_{\sigma(1)}W_{\sigma(2)}W_{\sigma(3)}}$.

Finally, it is clear from Proposition \ref{log:logwt}(b) that, in the
nontrivial logarithmic intertwining operator case, taking projections
to (generalized) weight subspaces is not enough to recover the
coefficients of $x^n(\log x)^k$ in ${\cal Y}(w_{(1)},x)w_{(2)}$ for
each $n\in\mathbb{C}$ and $k\in\mathbb{N}$.  However, by using
identities related to those in Lemma \ref{log:lemma} we have the
following result:
\begin{propo}
Let $W_1$, $W_2$, $W_3$ be generalized $V$-modules and ${\cal Y}$ a
logarithmic intertwining operator of type ${W_3\choose W_1\, W_2}$.
Let $w_{(1)}\in W_1$ and $w_{(2)}\in W_2$ be homogeneous of
(generalized) weights $n_1$ and $n_2$, respectively.
Then for any $n\in {\mathbb C}$ and any $r\in {\mathbb N}$,
${w_{(1)}}^{\cal Y}_{n;\,r}w_{(2)}$ can be written as a certain linear
combination of products of the component of weight $n_1+n_2-n-1$ of
\[
(L(0)-n_1-n_2+n+1)^l{\cal Y}((L(0)-n_1)^iw_{(1)},x)(L(0)-n_2)^jw_{(2)}
\]
for certain $i,j,l\in{\mathbb N}$ with monomials of the form
$x^{n+1}(\log x)^m$ for certain $m\in{\mathbb N}$.
\end{propo}

\setcounter{equation}{0}
\setcounter{rema}{0}

\section{Definition and constructions of $P(z)$-tensor product}
We now present our results on generalization of the $P(z)$-tensor
product construction developed in \cite{tensor1}, \cite{tensor2} and
\cite{tensor3}. Although many results here have the same statement as
in the ordinary case, one should note that the involvement of
logarithms makes the situation subtle, and their validity is based on
the results from the last section.

{\it Throughout this section and the remainder of this paper, we shall
assume the following, unless other assumptions are explicitly made:
$A$ is an abelian group and $\tilde{A}$ is an abelian group containing
$A$ as a subgroup; $V$ is a strongly $A$-graded conformal vertex
algebra; all $V$-modules and generalized $V$-modules considered are
strongly $\tilde{A}$-graded; and all intertwining operators and
logarithmic intertwining operators considered are grading-compatible.
Also, $z$ will be a fixed nonzero complex number.}

We first generalize the notion of $P(z)$-intertwining map from Section
4 of \cite{tensor1}:
\begin{defi}\label{im:imdef}{\rm
Let $(W_1,Y_1)$, $(W_2,Y_2)$ and $(W_3,Y_3)$ be generalized
$V$-modules.  A {\it $P(z)$-intertwining map of type ${W_3\choose
W_1\,W_2}$} is a linear map $I: W_1\otimes W_2 \to \overline{W}_3$
such that the following conditions are satisfied: the {\it grading
compatibility condition}: for $\beta, \gamma\in \tilde{A}$ and
$w_{(1)}\in W_{1}^{(\beta)}$, $w_{(2)}\in W_{2}^{(\gamma)}$,
\begin{equation}\label{grad-comp}
I(w_{(1)}\otimes w_{(2)})\in \overline{W_{3}^{(\beta+\gamma)}};
\end{equation}
the
{\em lower truncation condition:} for any elements $w_{(1)}\in W_1$,
$w_{(2)}\in W_2$, and any $n\in \mathbb{C}$,
\begin{equation}\label{im:ltc}
\pi_{n-m}I(w_{(1)}\otimes w_{(2)})=0\;\;\mbox{ for }\;m\in \mathbb{N}
\;\mbox{ sufficiently large}
\end{equation}
(which in fact follows from (\ref{grad-comp})); and the {\em Jacobi identity}:
\begin{eqnarray}\label{im:def}
\lefteqn{x_0^{-1}\delta\bigg(\frac{ x_1-z}{x_0}\bigg)
Y_3(v, x_1)I(w_{(1)}\otimes w_{(2)})}\nno\\
&&=z^{-1}\delta\bigg(\frac{x_1-x_0}{z}\bigg)
I(Y_1(v, x_0)w_{(1)}\otimes w_{(2)})\nno\\
&&\hspace{2em}+x_0^{-1}\delta\bigg(\frac{z-x_1}{-x_0}\bigg)
I(w_{(1)}\otimes Y_2(v, x_1)w_{(2)})
\end{eqnarray}
for $v\in V$, $w_{(1)}\in W_1$ and $w_{(2)}\in W_2$ (note that the
left-hand side of (\ref{im:def}) is meaningful because any infinite
linear combination of $v_n$ ($n\in\mathbb{Z}$) of the form
$\sum_{n<N}a_nv_n$ ($a_n\in \mathbb{C}$) acts on any $I(w_{(1)}\otimes
w_{(2)})$, due to (\ref{im:ltc})). The vector space of
$P(z)$-intertwining maps of type ${W_3}\choose {W_1W_2}$ is denoted by
${\cal M}[P(z)]^{W_3}_{W_1W_2}$, or simply by ${\cal
M}^{W_3}_{W_1W_2}$ if there is no ambiguity. }
\end{defi}

Following \cite{tensor1} we will choose the branch of $\log z$ so that
$0\leq {\rm Im }(\log z) < 2 \pi$. We will also use the notation
$l_p(z)=\log z+2\pi ip, \ p\in \mathbb{Z}$, from \cite{tensor1} for
arbitrary branches of the $\log$ function. For any expression $f(x)$
as in (\ref{log:f}) and any $\zeta\in \mathbb{C}$, whenever
\begin{equation}\label{log:fsub}
f(x)\lbar_{x^n=e^{n\zeta},\;(\log x)^m=\zeta^m,\;m,n\in\mathbb{C}}
\end{equation}
exists, in particular, when $f(x)={\cal Y}(w_{(1)},x)w_{(2)}$ $(\in
W_3[\log x]\{x\})$ for some $w_{(1)}\in W_1$, $w_{(2)}\in W_2$ and
some logarithmic intertwining operator ${\cal Y}$ of type ${W_3\choose
W_1\,W_2}$ ((\ref{log:fsub}) exists in this case because of
Proposition \ref{log:logwt}(b)), we will write (\ref{log:fsub}) simply
as $f(x)\lbar_{x=e^{\zeta}}$ or $f(e^\zeta)$, and call this
``substituting $e^\zeta$ for $x$ in $f(x)$'', despite the fact that
in general it depends on $\zeta$, not just $e^\zeta$ (see also
(\ref{log:subs})). In addition, fixing an integer $p$, we will sometimes
write
\[
f(x)\lbar_{x=z}\;\;\mbox{or}\;\;f(z)
\]
instead of $f(x)\lbar_{x=e^{l_p(z)}}$ or $f(e^{l_p(z)})$.

Fix an integer $p$. Let ${\cal Y}$ be a logarithmic intertwining
operator of type ${W_3\choose W_1\,W_2}$. Then we have a linear map
$I_{{\cal Y},p}: W_1\otimes W_2\to \overline{W}_3$ defined by
\begin{equation}\label{log:IYp}
I_{{\cal Y},p}(w_{(1)}\otimes w_{(2)})={\cal
Y}(w_{(1)},e^{l_p(z)})w_{(2)}
\end{equation}
for all $w_{(1)}\in W_1$ and $w_{(2)}\in W_2$ (the right-hand side of
(\ref{log:IYp}) is indeed an element of $\overline W_3$, again because
of Proposition \ref{log:logwt}(b)).  On the other hand, given a
$P(z)$-intertwining map $I$, we define a linear map ${\cal
Y}_{I,p}:W_1\otimes W_2\to W_3[\log x]\{x\}$ by
\begin{eqnarray}\label{recover}
\lefteqn{{\cal Y}_{I,p}(w_{(1)}, x)w_{(2)}=}\nno\\
&&y^{L(0)}x^{L(0)}I(y^{-L(0)}x^{-L(0)}w_{(1)}\otimes
y^{-L(0)}x^{-L(0)}w_{(2)})\lbar_{y=e^{-l_{p}(z)}}
\end{eqnarray}
for any $w_{(1)}\in W_1$ and $w_{(2)}\in W_2$ (this is well defined
and indeed maps to $W_3[\log x]\{x\}$, in view of (\ref{log:x^L(0)})
and (\ref{im:ltc})). Define the notation
${w_{(1)}}_{n;k}^{I,p}w_{(2)}\in W_3$ by
\[
{\cal Y}_{I,p}(w_{(1)}, x)w_{(2)}=\sum_{n\in\mathbb{C}}\sum_{k\in
\mathbb{N}} {w_{(1)}}_{n;k}^{I,p}w_{(2)} x^{-n-1}(\log x)^k.
\]
Observe that since the operator $x^{\pm L(0)}$ always increases the
power of $x$ in an expression homogeneous of (generalized) weight $n$
by $\pm n$, from (\ref{recover}) we have that
${w_{(1)}}_{n;k}^{I,p}w_{(2)}\in(W_3)_{[n_1+n_2-n-1]}$ for $w_{(1)}\in
(W_1)_{[n_1]}$ and $w_{(2)}\in (W_2)_{[n_2]}$.

Using results in the last section we have the following correspondence
between logarithmic intertwining operators and $P(z)$-intertwining
maps, generalizing Proposition 12.2 in \cite{tensor3}:
\begin{propo}\label{im:correspond}
For $p\in \mathbb{Z}$, the correspondence ${\cal Y}\mapsto I_{{\cal Y},
p}$ is a linear isomorphism from the space ${\cal V}^{W_3}_{W_1W_2}$
of logarithmic intertwining operators of type ${W_3\choose W_1\,W_2}$
to the space ${\cal M}^{W_3}_{W_1W_2}$ of $P(z)$-intertwining maps of
the same type. Its inverse map is given by $I\mapsto {\cal Y}_{I,p}$.
\end{propo}

The definition of $P(z)$-tensor product, in terms of the
$P(z)$-intertwining maps defined earlier, is the same as that in
\cite{tensor1}, except that here we take into account that the module
category under consideration can vary.  For this, recall the category
${\cal C}_1$ from the end of Section 2. We have:
\begin{defi}{\rm
For $W_1, W_2\in \ob{\cal C}_1$, a {\em $P(z)$-product of $W_1$ and
$W_2$} is an object $(W_3,Y_3)$ in ${\cal C}_1$ together with a
$P(z)$-intertwining map $I_3$ of type ${W_3\choose W_1\,W_2}$. We
denote it by $(W_3,Y_3;I_3)$ or simply by $(W_3,I_3)$. Let
$(W_4,Y_4;I_4)$ be another $P(z)$-product of $W_1$ and $W_2$. A {\em
morphism} from $(W_3,Y_3;I_3)$ to $(W_4,Y_4;I_4)$ is a module map
$\eta$ from $W_3$ to $W_4$ such that $I_4=\bar\eta\circ I_3$, where
$\bar\eta$ is the natural map {}from $\overline W_3$ to $\overline
W_4$ uniquely extending $\eta$.}
\end{defi}

\begin{defi}{\rm
Let ${\cal C}$ be a full subcategory of ${\cal C}_1$. For $W_1, W_2\in
\ob{\cal C}$, a {\em $P(z)$-tensor product of $W_1$ and $W_2$ in
${\cal C}$} is a $P(z)$-product $(W_0, Y_0; I_0)$ with $W_0\in{\rm
ob\,}{\cal C}$ such that for any $P(z)$-product $(W,Y;I)$ with
$W\in{\rm ob\,}{\cal C}$, there is a unique morphism from $(W_0, Y_0;
I_0)$ to $(W,Y;I)$.  Clearly, a $P(z)$-tensor product of $W_1$ and
$W_2$ in ${\cal C}$, if it exists, is unique up to a unique
isomorphism. In this case we will denote it as $(W_1\boxtimes_{P(z)}
W_2, Y_{P(z)}; \boxtimes_{P(z)})$ and call the object
$(W_1\boxtimes_{P(z)} W_2, Y_{P(z)})$ the {\em $P(z)$-tensor product
module of $W_1$ and $W_2$ in ${\cal C}$}.  We will omit the phrase
``in ${\cal C}$'' if the category ${\cal C}$ under consideration is
understood in context.  }
\end{defi}

\begin{rema}\label{anyz_1}{\rm
It is easy to show that in this setting, if $W_1\boxtimes_{P(z)} W_2$
exists, then $W_1\boxtimes_{P(z_1)} W_2$ exists {\it for any} $z_1 \in
\mathbb{C}^{\times}$.}
\end{rema}

Proposition \ref{im:correspond} shows that Proposition 12.3 in
\cite{tensor3}, which relates module maps from a $P(z)$-tensor product
module with $P(z)$-intertwining maps and intertwining operators, still
holds, for logarithmic intertwining operators in the present case.

Let $(W_1\boxtimes_{P(z)} W_2, Y_{P(z)}; \boxtimes_{P(z)})$ be the
$P(z)$-tensor product of $W_1$ and $W_2$.  We will sometimes denote
the image of the $P(z)$-intertwining map
\[
w_{(1)} \otimes w_{(2)}\mapsto\boxtimes_{P(z)}(w_{(1)} \otimes
w_{(2)})=\boxtimes_{P(z)}(w_{(1)}, z)w_{(2)}
\]
simply by $w_{(1)}\boxtimes_{P(z)} w_{(2)}$.

\begin{rema}\label{bifunctor}{\rm
Suppose that ${\cal C}$ is a full subcategory of ${\cal C}_1$ such
that for all $W_1, W_2 \in{\rm ob\,}{\cal C}$, the $P(z)$-tensor
product of $W_1$ and $W_2$ exists in ${\cal C}$.  Then the
$P(z)$-tensor product can be viewed as a (bi)functor from ${\cal C}
\times {\cal C}$ to ${\cal C}$ in a natural way.}
\end{rema}

As in the ordinary case, we have:
\begin{propo}
The module $W_1\boxtimes_{P(z)}W_2$ (if it exists) is spanned (as a
vector space) by the (generalized-)weight components of the elements
of $\overline{W_1\boxtimes_{P(z)}W_2}$ of the form
$w_{(1)}\boxtimes_{P(z)} w_{(2)}$, for all $w_{(1)}\in W_1$ and
$w_{(2)}\in W_2$.
\end{propo}

{\it Let ${\cal C}$ be some full subcategory of ${\cal C}_1$ closed
under the contragredient functor $(\cdot)'$ and under taking finite
direct sums, submodules and quotients, and containing $V$ as an
object.}

We now use the ``double dual'' approach, as in \cite{tensor1},
\cite{tensor2} and \cite{tensor3}, to construct the $P(z)$-tensor
product in ${\cal C}$, when it exists.  First we study a certain
adjoint of a $P(z)$-intertwining map.  Recall from \cite{FLM} that
$\iota_{+}$ is the map that expands a formal rational function (in
$t$, below) in the direction of nonnegative powers of the formal
variable.  Also recall from \cite{tensor1} the notation
\[
Y_{t}(v,x) = \sum_{n\in \mathbb{Z}} (v \otimes t^n)x^{-n-1}
\]
for $v\in V$.
\begin{defi}\label{tauPz}
{\rm Define the linear action $\tau_{P(z)}$ of
\[
V \otimes \iota_{+}\mathbb{C}[t,t^{- 1}, (z^{-1}-t)^{-1}]
\]
on $(W_1 \otimes W_2)^*$ by (cf. (5.1) in \cite{tensor1})
\begin{eqnarray*}
\lefteqn{\bigg(\tau_{P(z)}
\bigg(x_0^{-1}\delta\bigg(\frac{x^{-1}_1-z}{x_0}\bigg) Y_{t}(v,
x_1)\bigg)\lambda\bigg)(w_{(1)}\otimes w_{(2)})}\\
&&=z^{-1}\delta\bigg(\frac{x^{-1}_1-x_0}{z}\bigg)
\lambda(Y_1(e^{x_1L(1)}(-x_1^{-2})^{L(0)}v, x_0)w_{(1)}\otimes w_{(2)})\\
&&\quad +x^{-1}_0\delta\bigg(\frac{z-x^{-1}_1}{-x_0}\bigg)
\lambda(w_{(1)}\otimes Y_2^{o}(v, x_1)w_{(2)})
\end{eqnarray*}
for $v\in V$, $\lambda\in (W_1\otimes W_2)^{*}$, $w_{(1)}\in W_1$,
$w_{(2)}\in W_2$. (Note that this is indeed well defined and the
right-hand side is in fact a finite sum in view of the lower
truncation condition for vertex operators.)  Denote by $Y'_{P(z)}$ the
action of $V\otimes\mathbb{C}[t,t^{-1}]$ on $(W_1\otimes W_2)^*$ thus
defined: $Y'_{P(z)}(v,x)=\tau_{P(z)}(Y_t(v,x))$, that is,
\begin{eqnarray*}
\lefteqn{(Y'_{P(z)}(v,x_1)\lambda)(w_{(1)}\otimes w_{(2)})=
\lambda(w_{(1)}\otimes Y_2^{o}(v, x_1)w_{(2)})}\\
&&\quad+\res_{x_0}z^{-1}\delta\bigg(\frac{x^{-1}_1-x_0}{z}\bigg)
\lambda(Y_1(e^{x_1L(1)}(-x_1^{-2})^{L(0)}v, x_0)w_{(1)}\otimes
w_{(2)}).
\end{eqnarray*}
}
\end{defi}

\begin{rema}\label{I-intw}{\rm
Using the action $\tau_{P(z)}$, (\ref{im:def}) can be equivalently
written as
\[
\left(x_0^{-1}\delta\bigg(\frac{x^{-1}_1-z}{x_0}\bigg)Y'(v,x_1)
w'_{(3)}\right)\circ I=\tau_{P(z)}\bigg(x_0^{-1}\delta\bigg(\frac
{x^{-1}_1-z}{x_0}\bigg)Y_t(v, x_1)\bigg)(w'_{(3)}\circ I).
\]
for $w'_{(3)}\in W'_3$.}
\end{rema}

Using $\tau_{P(z)}$ we find that Proposition 13.3 in \cite{tensor3}
and the corresponding commutator formula still hold.

Write
\[
Y'_{P(z)}(\omega,x)=\sum_{n\in \mathbb{Z}} L'_{P(z)}(n)x^{-n-2}.
\]
Then we have that the coefficient operators of $Y'_{P(z)}(\omega, x)$
satisfy the Virasoro algebra commutator relations, that is,
\[
[L'_{P(z)}(m), L'_{P(z)}(n)]
=(m-n)L'_{P(z)}(m+n)+{\displaystyle\frac1{12}}(m^3-m)\delta_{m+n,0}c.
\]

Another notion, corresponding to the lower truncation condition for
$P(z)$-intertwining maps, is needed:
\begin{defi}{\rm
A map $J\in \hom(W'_3, (W_1\otimes W_2)^{*})$ is said to be {\em
grading restricted} if $J$ respects the $\tilde A$-gradings (the
notion of homogeneity of an element of $(W_1\otimes W_2)^*$ with
respect to $\tilde A$ is defined in the obvious way using the usual
tensor product $\tilde A$-grading of $W_1\otimes W_2$) and if given
any $n\in \mathbb{C}$, $w_{(1)}\in W_1$ and $w_{(2)}\in W_2$,
\[
J(({W'_3})_{[n-m]})(w_{(1)}\otimes w_{(2)})=0\;\;\mbox{ for }\;m\in
\mathbb{N}\;\mbox{ sufficiently large.}
\]
}
\end{defi}

Clearly, $J\in \hom(W'_3, (W_1\otimes W_2)^{*})$ is grading restricted
if and only if the map $I\in \hom(W_1\otimes W_2, \overline{W}_3)$
defined by
\[
\langle w'_{(3)}, I(w_{(1)}\otimes w_{(2)})\rangle=
\langle J(w'_{(3)}), w_{(1)}\otimes w_{(2)}\rangle
\]
satisfies the grading-compatibility condition (\ref{grad-comp}) and
the lower truncation condition (\ref{im:ltc}).

{}From this and Remark \ref{I-intw} we have the following result
generalizing Proposition 13.1 in \cite{tensor3}:
\begin{propo}\label{pz}
Let $W_1$, $W_2$ and $W_3$ be objects in ${\cal C}_1$. Then under the
natural map
\begin{equation}\label{nat}
\hom(W_1\otimes W_2, \overline{W}_3)\to \hom(W'_3, (W_1\otimes
W_2)^{*}),
\end{equation}
the $P(z)$-intertwining maps of type ${W_3}\choose {W_1W_2}$
correspond exactly to the grading restricted maps in $\hom(W'_3,
(W_1\otimes W_2)^{*})$ intertwining the actions of
\[
V \otimes \iota_{+}\mathbb{C}[t,t^{- 1},(z^{-1}-t)^{-1}]
\]
on $W'_3$ and $(W_1\otimes W_2)^{*}$.
\end{propo}

\begin{rema}{\rm
Note that in contrast with the case in \cite{tensor3}, the natural map
(\ref{nat}) is in general only an injection, rather than an
isomorphism.  This is because $\overline{W}_3$ is the completion of
$W_3$ with respect to the $\mathbb{C}$-grading and not the $\tilde
A$-grading.}
\end{rema}

For $W_1, W_2\in \ob{\cal C}_1$, a $P(z)$-product $(W,I)$ of $W_1$ and
$W_2$ with $W\in\ob{\cal C}_1$ can now be reformulated as an object
$(W,Y)$ in ${\cal C}_1$ equipped with a linear map $I$ from
$W_1\otimes W_2$ to $\overline W$ such that the corresponding linear
map
\begin{eqnarray*}
I': W'&\to&(W_1\otimes W_2)^*,\\
w'&\mapsto&w'\circ I
\end{eqnarray*}
is grading restricted and satisfies the intertwining conditions in
Proposition \ref{pz}.

Let $W_1, W_2 \in \ob{\cal C}$.  Define $W_1\hboxtr_{P(z)} W_2$ to be
the sum (or union) of all $I'(W')\subset(W_1\otimes W_2)^*$, where
$(W,I)$ is a $P(z)$-product of $W_1$ and $W_2$ with $W\in\ob{\cal
C}$. Then from the constructions and a proof analogous to that in the
ordinary case we have the following generalization of Proposition 13.7
in \cite{tensor3}:
\begin{propo}
If $(W_1\hboxtr_{P(z)} W_2, Y'_{P(z)})$ is an object of ${\cal C}$,
then denoting by $(W_1\boxtimes_{P(z)} W_2, Y_{P(z)})$ its
contragredient module, we have that the $P(z)$-tensor product in
${\cal C}$ exists and is $(W_1\boxtimes_{P(z)} W_2, Y_{P(z)}; i')$,
where $i$ is the natural inclusion from $W_1\hboxtr_{P(z)} W_2$ to
$(W_1\otimes W_2)^*$.  Conversely, if the $P(z)$-tensor product of
$W_1$ and $W_2$ in ${\cal C}$ exists, then $(W_1\hboxtr_{P(z)} W_2,
Y'_{P(z)})$ is an object in ${\cal C}$.
\end{propo}

Let $(W,I)$ be a $P(z)$-product of $W_1$ and $W_2$ and $w'\in W'$.
It is easy to see
that $I'(w')$ satisfies the following nontrivial and subtle conditions
on $\lambda \in (W_1\otimes W_2)^{*}$:
\begin{description}
\item{\bf The $P(z)$-compatibility condition}

(a) The {\em lower truncation condition}: For all $v\in V$, the formal
Laurent series $Y'_{P(z)}(v, x)\lambda$ involves only finitely many
negative powers of $x$.

(b) The following formula holds:
\begin{eqnarray}\label{cpb}
\lefteqn{\tau_{P(z)}\bigg(x_0^{-1}\delta\bigg(\frac{x^{-1}_1-z}{x_0}
\bigg)
Y_{t}(v, x_1)\bigg)\lambda}\nno\\
&&=x_0^{-1}\delta\bigg(\frac{x^{-1}_1-z}{x_0}\bigg)
Y'_{P(z)}(v, x_1)\lambda  \;\;\mbox{ for all }\;v\in V.
\end{eqnarray}
(Note that the two sides of (\ref{cpb}) are not {\it a priori} equal
for general $\lambda\in (W_1\otimes W_2)^{*}$.)
\end{description}

Furthermore, since $I'$ in particular intertwines the two actions of
$V\otimes\mathbb{C}[t, t^{-1}]$, $I'(W')$ is a generalized $V$-module.
Therefore, for each $w'\in
W'$, $I'(w')$ also satisfies the following condition on an
element $\lambda \in (W_1\otimes W_2)^*$:
\begin{description}
\item{\bf The $P(z)$-local grading restriction condition}

(a) The {\em grading condition}: $\lambda$ is a (finite) sum of
generalized eigenvectors in $(W_1\otimes W_2)^*$ for the operator
$L'_{P(z)}(0)$ that are also homogeneous with respect to $\tilde A$.

(b) The smallest $\tilde A$-graded subspace $W_\lambda$ of $(W_1\otimes W_2)^*$
containing $\lambda$ and stable under the component operators
$\tau_{P(z)}(v\otimes t^m)$ of the operators $Y'_{P(z)}(v,x)$ for
$v\in V$, $m\in \mathbb{Z}$ has the properties
\begin{eqnarray*}
&\dim(W_\lambda)^{(\beta)}_{[n]}<\infty ,&\\
&(W_\lambda)^{(\beta)}_{[n+k]}=0\;\;\mbox{ for }\;k\in \mathbb{Z}
\;\mbox{ sufficiently negative}&
\end{eqnarray*}
for any $n\in \mathbb{C}$ and $\beta\in \tilde A$, where the subscripts
denote the $\mathbb{C}$-grading by $L'_{P(z)}(0)$-eigenvalues and the
superscripts denote the natural $\tilde A$-grading.
\end{description}

The next two results are crucial to the theory, as the corresponding
results were in
\cite{tensor1}--\cite{tensor3}. Using the results on logarithmic
operators in the last section we are able to prove, by methods
analogous to, but more subtle than, the methods in
\cite{tensor1}--\cite{tensor3}:
\begin{theo}
Let $\lambda$ be an element of $(W_{1}\otimes W_{2})^{*}$ satisfying
the $P(z)$-compatibility condition. Then when acting on $\lambda$, the Jacobi
identity for $Y'_{P(z)}$ holds, that is,
\begin{eqnarray*}
\lefteqn{x_{0}^{-1}\delta
\left({\displaystyle\frac{x_{1}-x_{2}}{x_{0}}}\right)Y'_{P(z)}(u, x_{1})
Y'_{P(z)}(v, x_{2})\lambda}\\
&&\hspace{2ex}-x_{0}^{-1} \delta
\left({\displaystyle\frac{x_{2}-x_{1}}{-x_{0}}}\right)Y'_{P(z)}(v, x_{2})
Y'_{P(z)}(u, x_{1})\lambda\\
&&=x_{2}^{-1}\delta
\left({\displaystyle\frac{x_{1}-x_{0}}{x_{2}}}\right)Y'_{P(z)}(Y(u,
x_{0})v, x_{2})\lambda
\end{eqnarray*}
for $u, v\in V$.
\end{theo}

\begin{propo}
The subspace consisting of the elements of $(W_{1}\otimes W_{2})^{*}$
satisfying the $P(z)$-compatibility condition is stable under the operators
$\tau_{P(z)}(v\otimes t^{n})$ for $v\in V$ and $n\in \mathbb{Z}$, and
similarly for the subspace consisting of the elements satisfying the
$P(z)$-local grading-restriction condition.
\end{propo}

As a consequence of these results, we have that an element of
$(W_1\otimes W_2)^*$ satisfying the $P(z)$-compatibility condition and
the $P(z)$-local grading restriction condition generates a (strongly
graded) generalized $V$-module under the operators
$\tau_{P(z)}(v\otimes t^{n})$ for $v\in V$ and $n\in \mathbb{Z}$.

{}From all these results we can use the strategy of
\cite{tensor1}--\cite{tensor3} to obtain the following alternative
description of $W_1\hboxtr_{P(z)} W_2$:

\begin{propo}\label{characterizationofbackslash}
Suppose that for every element $\lambda$ of $(W_1\otimes W_2)^*$
satisfying the $P(z)$-compatibility condition and the $P(z)$-local
grading restriction condition, the generalized module $W_\lambda$ generated by
$\lambda$ under the operators $\tau_{P(z)}(v\otimes t^{n})$ for $v\in
V$ and $n\in \mathbb{Z}$ lies in ${\cal C}$ (this of course holds in
particular if ${\cal C}={\cal C}_1$). Then the subspace of
$(W_1\otimes W_2)^*$ consisting of all such elements is equal to
$W_1\hboxtr_{P(z)}W_2$.
\end{propo}

In the earlier work \cite{tensor1}--\cite{tensor3} of the first two
authors, another type of tensor product, the $Q(z)$-tensor product,
was studied and developed in addition to the $P(z)$-tensor product, in
\cite{tensor1} and \cite{tensor2}.  Then the $P(z)$-tensor product was
studied systematically in \cite{tensor3}, and many proofs for the
$P(z)$ case were obtained by the use of the results established for
the $Q(z)$ case.  Below we give the definition of a
$Q(z)$-intertwining map and its relations to $P(z)$-intertwining maps,
in the present more general context.

\begin{defi}{\rm
Let $(W_1,Y_1)$, $(W_2,Y_2)$ and $(W_3,Y_3)$ be generalized
$V$-modules.  A {\it $Q(z)$-intertwining map of type ${W_3\choose
W_1\,W_2}$} is a linear map $I: W_1\otimes W_2 \to \overline{W}_3$
such that the following conditions are satisfied: the {\it grading
compatibility condition}: for $\beta, \gamma\in \tilde{A}$ and
$w_{(1)}\in W_{1}^{(\beta)}$, $w_{(2)}\in W_{2}^{(\gamma)}$,
\begin{equation}
I(w_{(1)}\otimes w_{(2)})\in \overline{W_{3}^{(\beta+\gamma)}};
\end{equation}
the {\em lower truncation condition:} for any elements $w_{(1)}\in
W_1$, $w_{(2)}\in W_2$, and any $n\in \mathbb{C}$,
\begin{equation}\label{imq:ltc}
\pi_{n-m}I(w_{(1)}\otimes w_{(2)})=0\;\;\mbox{ for }\;m\in
\mathbb{N}\;\mbox{ sufficiently large;}
\end{equation}
and the {\em Jacobi identity}:
\begin{eqnarray}\label{imq:def}
\lefteqn{z^{-1}\delta\left(\frac{x_1-x_0}{z}\right)
Y^o_3(v, x_0)I(w_{(1)}\otimes w_{(2)})}\nonumber\\
&&=x_0^{-1}\delta\left(\frac{x_1-z}{x_0}\right)
I(Y_1^{o}(v, x_1)w_{(1)}\otimes w_{(2)})\nonumber\\
&&\hspace{2em}-x_0^{-1}\delta\left(\frac{z-x_1}{-x_0}\right)
I(w_{(1)}\otimes Y_2(v, x_1)w_{(2)})
\end{eqnarray}
for $v\in V$, $w_{(1)}\in W_1$ and $w_{(2)}\in W_2$ (recall
(\ref{yo}), and note that the left-hand side of (\ref{imq:def}) is
meaningful because any infinite linear combination of $v_n$ of the
form $\sum_{n<N}a_nv_n$ ($a_n\in \mathbb{C}$) acts on any
$I(w_{(1)}\otimes w_{(2)})$, due to (\ref{imq:ltc})). The vector space
of $Q(z)$-intertwining maps of type ${W_3\choose W_1\,W_2}$ is denoted
by ${\cal M}[Q(z)]^{W_3}_{W_1W_2}$, or simply by ${\cal
M}^{W_3}_{W_1W_2}$ if there is no ambiguity. }
\end{defi}

We have the following result relating $P(z)$- and $Q(z)$-intertwining
maps, essentially due to the fact that $Q(z)$ represents the Riemann
sphere $\mathbb{C}\cup \{ \infty \}$ with ordered punctures $z$,
$\infty$, $0$, {\it with $z$ the negatively oriented puncture}, and
with standard local coordinates vanishing at these punctures, while
$P(z)$ represents the Riemann sphere with ordered punctures $\infty$,
$z$, $0$, {\it with $\infty$ the negatively oriented puncture}, and
with standard local coordinates vanishing at the punctures:

\begin{propo}
Let $W_1$, $W_2$ and $W_3$ be generalized $V$-modules. Let $I:
W_1\otimes W_2\to \overline{W}_3$ and $J: W'_3\otimes W_2\to
\overline{W'_1}$ be linear maps related to each other by
\[
\langle w_{(1)},J(w'_{(3)}\otimes w_{(2)})\rangle=
\langle w'_{(3)},I(w_{(1)}\otimes w_{(2)})\rangle
\]
for any $w_{(1)}\in W_1$, $w_{(2)}\in W_2$ and $w'_{(3)}\in
W'_3$. Then $I$ is a $Q(z)$-intertwining map of type ${W_3\choose
W_1\,W_2}$ if and only if $J$ is a $P(z)$-intertwining map of type
${W'_1\choose W'_3\,W_2}$.
\end{propo}

We remark that all the material in the earlier part of this section
can be carried over to the $Q(z)$ case.

\setcounter{equation}{0}
\setcounter{rema}{0}

\section{The associativity isomorphism}
Having developed the concept of logarithmic intertwining operator and
constructed the $P(z)$-tensor product, when it exists, we now proceed
to the construction of the associativity isomorphism for this type of
tensor product.  We are able to show that the material in
\cite{tensor4} can be carried over, or adapted, to the setup in this
paper. We first investigate the necessary conditions for the existence
of the associativity isomorphism, and then we give the construction
under these conditions. Again we refer the reader to \cite{HLZ} for
further details and proofs.

We adopt the definition (incorporating convergence) of existence of
products (respectively, iterates) of two intertwining maps given in
(14.1) (respectively, (14.2)) in \cite{tensor4}. We have the following
result generalizing that in \cite{tensor4}, which can be proved by a
logarithmic analogue of the arguments in \cite{tensor4}.  One needs to
use the transformations $\Omega_r$ developed in Section 3, which
involve $\log x$ in a subtle way.
\begin{propo}\label{convergence}
The following two conditions are equivalent:
\begin{enumerate}
\item For any objects $W_1$, $W_2$, $W_3$, $W_4$ and $M_1$ in ${\cal
C}$, any nonzero complex numbers $z_1$ and $z_2$ satisfying
$|z_1|>|z_2|>0$, any $P(z_1)$-intertwining map $I_1$ of type
${W_4}\choose {W_1M_1}$ and any $P(z_2)$-intertwining map $I_2$ of
type ${M_1}\choose {W_2W_3}$, the product $\gamma (I_1;1_{W_1},I_2)$
of $I_1$ and $I_2$ exists.
\item For any objects $W_1$, $W_2$, $W_3$, $W_4$ and $M_2$ in ${\cal
C}$, any nonzero complex numbers $z_0$ and $z_2$ satisfying
$|z_2|>|z_0|>0$, any $P(z_2)$-intertwining map $I^1$ of type
${W_4}\choose {M_2W_3}$ and any $P(z_0)$-intertwining map $I^2$ of
type ${M_2}\choose {W_1W_2}$, the iterate $\gamma(I^1;I^2,1_{W_3})$ of
$I^1$ and $I^2$ exists.
\end{enumerate}
\end{propo}

We call either property in Proposition \ref{convergence} the {\it
convergence condition in the category ${\cal C}$}. If
the convergence condition in the category ${\cal C}$ holds, then
as images of the product or iterate of certain intertwining maps,
$w_{(1)}\boxtimes_{P(z_1)}(w_{(2)}\boxtimes_{P(z_2)}w_{(3)})$ and
$(w_{(1)}\boxtimes_{P(z_1-z_{2})}w_{(2)})\boxtimes_{P(z_2)}w_{(3)}$
for $w_{(1)}\in W_{1}$, $w_{(2)}\in W_{2}$ and $w_{(3)}\in W_{3}$
exist as elements of
\[
\overline{W_{1}\boxtimes_{P(z_1)}(W_{2}\boxtimes_{P(z_2)}W_{3})}
\;\;\mbox{and}\;\;
\overline{(W_{1}\boxtimes_{P(z_1-z_{2})}W_{2})\boxtimes_{P(z_2)}W_{3}},
\]
respectively, when these tensor products exist.

We will use the following concept concerning unique expansion of an
analytic function in terms of powers of $z$ and $\log z$: Here we will
call a subset ${\cal S}$ of $\mathbb{C}\times\mathbb{C}$ a {\em unique
expansion set} if the absolute convergence to $0$ on some nonempty
open subset of $\mathbb{C}$ of any series
$\sum_{(\alpha,\beta)\in{\cal S}} a_{\alpha,\beta}z^\alpha(\log
z)^\beta$ implies that $a_{\alpha,\beta}=0$ for all
$(\alpha,\beta)\in{\cal S}$.  It is easy to show that
$\mathbb{Z}\times\{ 0,1,\dots,N\}$ is a unique expansion set for any
$N\in\mathbb{N}$.  More generally, $D\times\{0,1,\dots,N\}$ is a
unique expansion set for any discrete subset $D$ of $\mathbb{R}$.  On
the other hand, it is known that $\mathbb{C}\times\{0\}$ is {\em not}
a unique expansion set\footnote{We thank A.~Eremenko for informing us
of this result.}.

{\it For the rest of this section, we assume for convenience that for
any object in ${\cal C}$, the (generalized) weights form a discrete
set of real numbers.}  This assumption implies that for any objects
$W_1$, $W_2$ and $W_3$ in ${\cal C}$, all possible powers (with
nonzero coefficients) of $x$ and $\log x$ in the image of any element
of $W_{1}\otimes W_{2}$ under any logarithmic intertwining operator of
type ${W_3\choose W_1\;W_2}$ form a unique expansion set.

Under this assumption, by using (\ref{log:ck1}) and
(\ref{log:L(-1)dev}) we can prove the following result:
\begin{propo}
Assume that the convergence condition in ${\cal C}$ holds.  Let $z_1$,
$z_2$ be two fixed nonzero complex numbers satisfying $|z_1|>|z_2|>0$,
and let $I_1$ and $I_2$ be $P(z_1)$- and $P(z_2)$-intertwining maps of
types ${W_4 \choose W_1\, M_1}$ and ${M_1 \choose W_2\, W_3}$,
respectively. Let $w_{(1)}\in W_1$, $w_{(3)}\in W_3$ and $w'_{(4)}\in
W'_4$. Suppose that for any $w_{(2)}\in W_2$,
\[
\langle w'_{(4)}, I_1(w_{(1)},z_1)I_2(w_{(2)},z_2)w_{(3)}\rangle=0.
\]
Then
\[
\langle w'_{(4)}, \pi_pI_1({w_{(1)}},z_1)\pi_q I_2({w_{(2)}},z_2)
w_{(3)}\rangle=0
\]
(recall (\ref{pin})) for all $p, q\in \mathbb{C}$ and $w_{(2)}\in W_2$.
\end{propo}

As a consequence, we have:
\begin{corol}
Suppose that the $P(z_2)$-tensor product of $W_2$ and $W_3$ and the
$P(z_1)$-tensor product of $W_1$ and $W_2 \boxtimes_{P(z_2)} W_3$ both
exist. Then the space $W_1\boxtimes_{P(z_1)} (W_2\boxtimes_{P(z_2)}
W_3)$ is spanned as a vector space by the weight components of all
elements of the form
$w_{(1)}\boxtimes_{P(z_1)}(w_{(2)}\boxtimes_{P(z_2)}w_{(3)})$, where
$w_{(1)}\in W_1$, $w_{(2)}\in W_2$ and $w_{(3)}\in W_3$.
\end{corol}

Analogously, we have the corresponding proposition and corollary
concerning iterates of intertwining maps as well, which we will not
give here.

\begin{rema}{\rm
One can generalize these propositions to the case of products and
iterates of more than two intertwining maps. The spanning property in
the case of four generalized modules will be used in proving the
``coherence theorem'' for the constructed ``vertex tensor category.''}
\end{rema}

{\it Now suppose that the convergence condition in ${\cal C}$ holds.}
We first give the following definition, which will be important in
studying compositions of intertwining maps:
\begin{defi}
{\rm Let $z_0, z_1, z_2\in \mathbb{C}^{\times }$ with $z_0=z_1-z_2$
(so that in particular $z_1\neq z_2$, $z_0\neq z_1$ and $z_0\neq
-z_2$).  Let $(W_1,Y_1)$, $(W_2,Y_2)$, $(W_3,Y_3)$ and ($W_4,Y_4)$ be
generalized $V$-modules.  A {\it $P(z_1,z_2)$-intertwining map} is a
linear map $F:\, W_1\otimes W_2\otimes W_3\to \overline{W}_4$ such
that the following conditions are satisfied: the {\it grading
compatibility condition}: for $\beta_i \in \tilde{A}$ and $w_{(i)}\in
W_{i}^{(\beta_i)}$ for $i=1,2,3$,
\begin{equation}\label{grcompcondn}
F(w_{(1)}\otimes w_{(2)}\otimes w_{(3)})\in \overline{W_{4}^{(\beta_1
+ \beta_2 + \beta_3)}};
\end{equation}
the {\em lower truncation condition}: for any $w_{(1)}\in W_1$,
$w_{(2)}\in W_2$, $w_{(3)}\in W_3$ and $n\in \mathbb{C}$,
\begin{equation}\label{zz:ltc}
\pi_{n-m}F(w_{(1)}\otimes w_{(2)}\otimes w_{(3)})=0\;\;\mbox{ for
}\;m\in \mathbb{N}\;\mbox{ sufficiently large;}
\end{equation}
and the {\em composite Jacobi identity}:
\begin{eqnarray}\label{zz:Y}
\lefteqn{\dlt{x_1}{x_0}{-z_1}\dlt{x_2}{x_0}{-z_2}Y_4(v,x_0)F(w_{(1)}
\otimes w_{(2)}\otimes w_{(3)})}\nno \\
&& =\dlt{x_0}{z_1}{+x_1}\dlt{x_2}{z_0}{+x_1}F(Y_1(v,x_1)w_{(1)}
\otimes w_{(2)}\otimes w_{(3)})\nno \\
&& +\dlt{x_0}{z_2}{+x_2}\dlt{x_1}{-z_0}{+x_2}F(w_{(1)}\otimes
Y_2(v,x_2)w_{(2)}\otimes w_{(3)})\nno \\
&& +\dlt{x_1}{-z_1}{+x_0}\dlt{x_2}{-z_2}{+x_0}F(w_{(1)}\otimes
w_{(2)}\otimes Y_3(v,x_0)w_{(3)})\nno \\
\end{eqnarray}
for $v\in V$, $w_{(1)}\in W_1$, $w_{(2)}\in W_2$ and $w_{(3)}\in W_3$
(note that the left-hand side of (\ref{zz:Y}) is meaningful because
any infinite linear combination of the $v_n$ of the form
$\sum_{n<N}a_nv_n$ ($a_n\in \mathbb{C}$) acts on any $F(w_{(1)}\otimes
w_{(2)}\otimes w_{(3)})$, due to (\ref{zz:ltc})).  }
\end{defi}

As in \cite{tensor4}, one shows that in the setting of Proposition
\ref{convergence}, $\gamma (I_1;1_{W_1},I_2)$ is a
$P(z_1,z_2)$-intertwining map when $|z_1|>|z_2|>0$, and
$\gamma(I^1;I^2,1_{W_3})$ is a $P(z_2+z_0,z_2)$-intertwining map when
$|z_2|>|z_0|>0$.

We now define the following action (cf. Definition \ref{tauPz}):
\begin{defi}{\rm
Let $z_1, z_2\in \mathbb{C}^{\times}$, $z_1\neq z_2$. The linear action
$\tau_{P(z_1,z_2)}$ of
\begin{equation}\label{Vtensoriota}
V\otimes \iota _{+}\mathbb{C}[t,t^{-1},(z_1^{-1}-t)^{-1},(z_2^{-1}
-t)^{-1}]
\end{equation}
on $(W_1\otimes W_2\otimes W_3)^{*}$ is defined by
\begin{eqnarray*}
\lefteqn{\left(\tau_{P(z_1,z_2)}\left(\dlt{x_1}{x_0^{-1}}{-z_1}\dlt{x_2}
{x_0^{-1}}{-z_2}Y_{t}(v,x_0)\right)\lambda \right)(w_{(1)}\otimes
w_{(2)}\otimes w_{(3)})}\\
&=& \dlti{x_0}{z_1}{+x_1}\dlt{x_2}{z_0}{+x_1}\cdot\\
&& \qquad\qquad \lambda (Y_1((-x_0^{-2})^{L(0)}e^{-x_0^{-1}L(1)}v,x_1)
w_{(1)}\otimes w_{(2)}\otimes w_{(3)})\\
&& +\dlti{x_0}{z_2}{+x_2}\dlt{x_1}{-z_0}{+x_2}\cdot\\
&& \qquad\qquad \lambda (w_{(1)}\otimes
Y_2((-x_0^{-2})^{L(0)}e^{-x_0^{-1}L(1)}v,x_2)w_{(2)}\otimes w_{(3)})\\
&& +\dlt{x_1}{-z_1}{+x^{-1}_0}\dlt{x_2}{-z_2}{+x^{-1}_0}\lambda
(w_{(1)}\otimes w_{(2)}\otimes Y_3^o(v,x_0)w_{(3)})
\end{eqnarray*}
for $v\in V$, $\lambda \in (W_1\otimes W_2\otimes W_3)^{*}$,
$w_{(1)}\in W_1$, $w_{(2)}\in W_2$ and $w_{(3)}\in W_3$.  Also, denote
by $Y'_{P(z_1,z_2)}$ the action of $V\otimes\mathbb{C}[t,t^{-1}]$ on
$(W_1\otimes W_2\otimes W_3)^*$ thus defined, that is,
\[
Y'_{P(z_1,z_2)}(v,x)=\tau_{P(z_1,z_2)}(Y_t(v,x)).
\]
}
\end{defi}

\begin{rema}{\rm
In (14.18) and (14.20) of \cite{tensor4}, actions $\tau^{(1)}_
{P(z_1,z_2)}$ and $\tau^{(2)}_{P(z_1,z_2)}$ of (\ref{Vtensoriota}) on
$(W_1 \otimes W_2 \otimes W_3)^*$ were defined, when $|z_1|>|z_2|>0$
and $|z_2|>|z_1-z_2|>0$, respectively.  The action $\tau
_{P(z_1,z_2)}$ of (\ref{Vtensoriota}) on $(W_1 \otimes W_2 \otimes
W_3)^*$, defined on all of $\{z_1,z_2\in \mathbb{C}^{\times}| z_1\neq
z_2\}$, coincides with, and thus extends, these two actions.}
\end{rema}

\begin{rema}{\rm
Using the action $\tau_{P(z_1,z_2)}$, the equality (\ref{zz:Y}) can be
equivalently written as: For any $w'_{(4)} \in W'_4$,
\begin{eqnarray*}
\lefteqn{\left(\dlt{x_1}{x^{-1}_0}{-z_1}\dlt{x_2}{x^{-1}_0}{-z_2}
Y'_4(v,x_0)w'_{(4)}\right)\circ F}\\
&&=\tau_{P(z_1,z_2)}\left(\dlt{x_1}{x_0^{-1}}{-z_1}\dlt{x_2}
{x_0^{-1}}{-z_2}Y_{t}(v,x_0)\right)(w'_{(4)}\circ F).
\end{eqnarray*}
}
\end{rema}

We will write
\[
Y'_{P(z_1,z_2)}(\omega ,x)=\sum _{n\in\mathbb{Z}}L'_{P(z_1,z_2)}(n)
x^{-n-2}.
\]

As in the $P(z)$-intertwining map case, one more notion, corresponding
to the lower truncation condition for $P(z_1,z_2)$-intertwining maps,
is needed:
\begin{defi}{\rm
A map $G\in \hom(W'_4, (W_1\otimes W_2\otimes W_3)^{*})$ is said to be
{\em grading restricted} if $G$ respects the $\tilde{A}$-gradings and
if given any $n\in \mathbb{C}$, $w_{(1)}\in W_1$, $w_{(2)}\in W_2$ and
$w_{(3)}\in W_3$,
\[
G((W'_4)_{[n-m]})(w_{(1)}\otimes w_{(2)}\otimes w_{(3)})=0\;\;\mbox{
for }\;m\in \mathbb{N}\;\mbox{ sufficiently large.}
\]
}
\end{defi}

Clearly, if $G\in \hom(W'_4, (W_1\otimes W_2\otimes W_3)^{*})$ is
grading restricted, then the map $F\in \hom(W_1\otimes W_2\otimes W_3,
\overline{W}_4)$ defined by
\[
\langle w'_{(4)}, F(w_{(1)}\otimes w_{(2)}\otimes w_{(3)}\rangle=
\langle G(w'_{(4)}), w_{(1)}\otimes w_{(2)}\otimes w_{(3)}\rangle
\]
satisfies the grading-compatibility condition (\ref{grcompcondn}) and
the lower truncation condition (\ref{zz:ltc}).

{}From these considerations we obtain:
\begin{propo}\label{zzcor}
Let $z_1,z_2\in \mathbb{C}^{\times }$, $z_1\neq z_2$. Let $W_1$, $W_2$,
$W_3$ and $W_4$ be objects in ${\cal C}$. Then under the natural map
\[
\hom(W_1\otimes W_2\otimes W_3,\overline{W}_4)\to
\hom(W_4',(W_1\otimes W_2\otimes W_3)^{*}),
\]
the $P(z_1,z_2)$-intertwining maps correspond exactly to the grading
restricted maps in $\hom(W_4',(W_1\otimes W_2\otimes W_3)^{*})$
intertwining the actions of the space (\ref{Vtensoriota}) on $W_4'$
and $(W_1\otimes W_2\otimes W_3)^*$.
\end{propo}

Let $G\in \hom (W'_4, (W_1\otimes W_2\otimes W_3)^*)$ correspond to a
$P(z_1, z_2)$-intertwining map as in Proposition \ref{zzcor}.  Then
for any $w'_{(4)}\in W'_4$, $G(w'_{(4)})$ satisfies the following
conditions on an element $\lambda \in (W_1\otimes W_2\otimes W_3)^*$:
\begin{description}
\item{\bf The $P(z_1, z_2)$-compatibility condition}

(a) The \emph{$P(z_1, z_2)$-lower truncation condition}: For any $v\in
V$, the formal Laurent series $Y'_{P(z_1,z_2)}(v,x)\lambda$ involves
only finitely many negative powers of $x$.

(b) The following formula holds: for any $v\in V$,
\begin{eqnarray}\label{zz:cpb}
\lefteqn{\tau_{P(z_1,z_2)}\bigg(\dlt{x_1}{x_0^{-1}}{-z_1}\dlt{x_2}
{x_0^{-1}}{-z_2}Y_{t}(v,x_0)\bigg)\lambda}\nno\\
&&=\dlt{x_1}{x_0^{-1}}{-z_1}\dlt{x_2}{x_0^{-1}}{-z_2}Y'_{P(z_1,z_2)}
(v,x_0)\lambda.
\end{eqnarray}
(Note that the two sides of (\ref{zz:cpb}) are not {\it a priori}
equal for general $\lambda\in (W_1\otimes W_2\otimes W_3)^{*}$.)

\item{\bf The $P(z_1,z_2)$-local grading restriction condition}

(a) The {\em grading condition}: $\lambda$ is a (finite) sum of
generalized eigenvectors in $(W_1\otimes W_2\otimes W_3)^{*}$ for the
operator $L'_{P(z_1,z_2)}(0)$ that are also homogeneous with respect
to $\tilde A$.

(b) The smallest $\tilde A$-graded subspace $W_\lambda$ of $(W_1\otimes W_2\otimes
W_3)^{*}$ containing $\lambda$ and stable under the component
operators $\tau_{P(z_1,z_2)}(v\otimes t^{m})$ of the operators
$Y'_{P(z_1,z_2)}(v,x)$ for $v\in V$, $m\in \mathbb{Z}$ has the
properties
\begin{eqnarray}
&\dim(W_\lambda)^{(\beta)}_{[n]}<\infty,&\\
&(W_\lambda)^{(\beta)}_{[n+k]}=0\;\;\mbox{ for }\;k\in \mathbb{Z}
\;\mbox{ sufficiently negative}&
\end{eqnarray}
for any $n\in \mathbb{C}$ and $\beta\in \tilde A$, where the subscripts
denote the $\mathbb{C}$-grading by $L'_{P(z_1,z_2)}(0)$-eigenvalues and
the superscripts denote the natural $\tilde A$-grading.
\end{description}

Recall that we have assumed the convergence condition in the category
${\cal C}$.  We shall now study the condition for the product of two
suitable intertwining maps to be written as the iterate of some
suitable intertwining maps, and vice versa.

Let $z_1$, $z_2$ be distinct nonzero complex numbers, and let
$z_0=z_1-z_2$.  Let $W_!$, $W_2$, $W_3$, $W_4$, $M_1$ and $M_2$ be
objects in ${\cal C}$ and let $I_1$, $I_2$, $I^1$ and $I^2$ be
$P(z_1)$-, $P(z_2)$-, $P(z_2)$- and $P(z_0)$-intertwining maps of
types ${W_4 \choose W_1\, M_1}$, ${M_1 \choose W_2\, W_3}$, ${W_4
\choose M_2\, W_3}$ and ${M_2\choose W_1\, W_2}$, respectively. Then
by the assumption of the convergence condition in ${\cal C}$, when
$|z_1|>|z_2|>|z_0|>0$, both $\gamma(I_1; 1_{W_1}, I_2)$ and
$\gamma(I^1; I^2, 1_{W_3})$ exist and are $P(z_1,z_2)$-intertwining
maps. We shall now discuss when two such $P(z_1,z_2)$-intertwining
maps are in fact equal to each other.

As in \cite{tensor4}, we need the following:
\begin{defi}{\rm
Let $\lambda \in (W_1\otimes W_2\otimes W_3)^*$.  Define
$\mu^{(1)}_{\lambda, w_{(1)}}$ to be the linear functional
$\lambda(w_{(1)}\otimes\cdot)\in
(W_2\otimes W_3)^{*}$ for any $w_{(1)}\in W_1$ and
$\mu^{(2)}_{\lambda, w_{(3)}}$ to be $\lambda(\cdot\otimes w_{(3)})\in
(W_1\otimes W_2)^{*}$ for any $w_{(3)}\in W_3$. }
\end{defi}

We have that Lemma 14.3 in \cite{tensor4} still holds in the
generality of the present paper.

As in \cite{tensor4} the following two conditions on $\lambda\in
(W_1\otimes W_2\otimes W_3)^*$ will be needed:
\begin{description}
\item{\bf $P^{(1)}(z)$-local grading restriction condition}

(a) The {\em $P^{(1)}(z)$-grading condition}: For any $w_{(1)}\in
W_1$, the element $\mu^{(1)}_{\lambda, w_{(1)}}\in (W_2\otimes
W_3)^{*}$ is the limit, in the locally convex topology defined by the
pairing between $(W_2\otimes W_3)^{*}$ and $W_2\otimes W_3$, of an
absolutely convergent series of generalized eigenvectors in
$(W_2\otimes W_3)^{*}$ with respect to the operator $L'_{P(z)}(0)$
that are also homogeneous with respect to $\tilde A$.

(b) For any $w_{(1)}\in W_1$, let $W^{(1)}_{\lambda, w_{(1)}}$ be the
smallest subspace of $(W_2\otimes W_3)^{*}$ containing the terms in
the series in (a) and stable under the component operators
$\tau_{P(z)}(v\otimes t^{m})$ of the operators $Y'_{P(z)}(v, x)$ for
$v\in V$, $m\in \mathbb{Z}$.  Then $W^{(1)}_{\lambda, w_{(1)}}$ has the
properties
\begin{eqnarray*}
&\dim(W^{(1)}_{\lambda, w_{(1)}})^{(\beta)}_{[n]}<\infty,&\\
&(W^{(1)}_{\lambda, w_{(1)}})^{(\beta)}_{[n+k]}=0\;\;\mbox{ for
}\;k\in \mathbb{Z} \;\mbox{ sufficiently negative}&
\end{eqnarray*}
for any $n\in \mathbb{C}$ and $\beta\in \tilde A$, where the subscripts
denote the $\mathbb{C}$-grading by $L'_{P(z)}(0)$-eigenvalues and the
superscripts denote the natural $\tilde A$-grading.

\item{\bf $P^{(2)}(z)$-local grading restriction condition}

(a) The {\em $P^{(2)}(z)$-grading condition}: For any $w_{(3)}\in W_3$,
the element $\mu^{(2)}_{\lambda, w_{(3)}}\in (W_1\otimes W_2)^{*}$ is
the limit, in the locally convex topology defined by the
pairing between $(W_1\otimes W_2)^{*}$ and $W_1\otimes W_2$,
of an absolutely convergent series of generalized
eigenvectors in $(W_1\otimes W_2)^{*}$ with respect to the operator
$L'_{P(z)}(0)$ that are also homogeneous with respect to $\tilde A$.

(b) For any $w_{(3)}\in W_3$, let $W^{(2)}_{\lambda, w_{(3)}}$ be the
smallest subspace of $(W_1\otimes W_2)^{*}$ containing the terms in
the series in (a) and stable under the component operators
$\tau_{P(z)}(v\otimes t^{m})$ of the operators $Y'_{P(z)}(v, x)$ for
$v\in V$, $m\in \mathbb{Z}$. Then $W^{(2)}_{\lambda, w_{(3)}}$ has the
properties
\begin{eqnarray*}
&\dim(W^{(2)}_{\lambda, w_{(3)}})^{(\beta)}_{[n]}<\infty,&\\
&(W^{(2)}_{\lambda, w_{(3)}})^{(\beta)}_{[n+k]}=0\;\;\mbox{ for
}\;k\in \mathbb{Z} \;\mbox{ sufficiently negative}&
\end{eqnarray*}
for any $n\in \mathbb{C}$ and $\beta\in \tilde A$, where the subscripts
denote the $\mathbb{C}$-grading by $L'_{P(z)}(0)$-eigenvalues and the
superscripts denote the natural $\tilde A$-grading.
\end{description}

By generalizing the arguments in \cite{tensor4}, we can now prove the
following two results:
\begin{propo}\label{5.12}
Let $I_{1}$, $I_{2}$, $I^1$ and $I^2$ be $P(z_1)$-, $P(z_2)$-, $P(z_2)$- and
$P(z_0)$-intertwining maps of types ${W_4}\choose {W_1M_1}$,
${M_1}\choose {W_2W_3}$, ${W_4}\choose {M_2W_3}$ and ${M_2}\choose
{W_1W_2}$, respectively. Then for any $w'_{(4)}\in W'_4$,
$\gamma(I_1;
1_{W_1}, I_2)'(w'_{(4)})$ satisfies the $P^{(1)}(z_2)$-local grading
restriction condition  when $|z_1|>|z_2|>0$, and
$\gamma(I^1;I^2,1_{W_3})'(w'_{(4)})$ satisfies the
$P^{(2)}(z_0)$-local grading restriction condition when
$|z_2|>|z_0|>0$.
\end{propo}
\begin{propo}\label{associativityoflogintwops}
Let $W_1$, $W_2$, $W_3$, $W'_4$ and $M_1$ be objects in ${\cal C}$
and assume that $W_1\boxtimes_{P(z_0)} W_2$ exists in ${\cal C}$.  Let
$I_1$ and $I_2$ be $P(z_1)$- and $P(z_2)$-intertwining maps of types
${W_4 \choose W_1\, M_1}$ and ${M_1 \choose W_2\, W_3}$, respectively.
Suppose that $\gamma(I_1; 1_{W_1}, I_2)'(w'_{(4)})$ satisfies the
$P^{(2)}(z_0)$-local grading restriction condition for all
$w'_{(4)}\in W'_4$.  Then there is a $P(z_2)$-intertwining map of type
${W_4\choose W_1\boxtimes_{P(z_0)} W_2\,\,W_3}$ such that
\[
\langle w'_{(4)},I_1(w_{(1)},z_1)I_2(w_{(2)},z_2)w_{(3)}\rangle
=\langle w'_{(4)},I(w_{(1)}\boxtimes_{P(z_0)} w_{(2)},z_2)
w_{(3)}\rangle
\]
for any $w_{(1)}\in W_1$, $w_{(2)}\in W_2$, $w_{(3)}\in W_3$ and
$w'_{(4)}\in W'_4$.  Conversely, let $W_1$, $W_2$, $W_3$, $W'_4$ and
$M_2$ be objects in ${\cal C}$ and assume that $W_2\boxtimes_{P(z_2)}
W_3$ exists in ${\cal C}$.  Let $I^1$ and $I^2$ be $P(z_2)$- and
$P(z_0)$-intertwining maps of types ${W_4 \choose M_2\, W_3}$ and
${M_2 \choose W_1\, W_2}$, respectively.  Suppose that $\gamma(I^1;
I^2, 1_{W_3})'(w'_{(4)})$ satisfies the $P^{(1)}(z_2)$-local grading
restriction condition for all $w'_{(4)}\in W'_4$.  Then there is a
$P(z_1)$-intertwining map of type ${W_4\choose
W_1\,\,W_2\boxtimes_{P(z_2)} W_3}$ such that
\[
\langle w'_{(4)},I^1(I^2(w_{(1)},z_0)w_{(2)},z_2)w_{(3)}\rangle
=\langle w'_{(4)},I(w_{(1)},z_1)(w_{(2)}\boxtimes_{P(z_2)}w_{(3)})
\rangle
\]
for any $w_{(1)}\in W_1$, $w_{(2)}\in W_2$, $w_{(3)}\in W_3$ and
$w'_{(4)}\in W'_4$.
\end{propo}

Proposition 14.7 and Theorem 14.8 in \cite{tensor4} also hold in the
present generality.  In particular, Proposition 14.7 in \cite{tensor4}
in our setting states:
\begin{propo}\label{expansion}
The following two conditions are equivalent:
\begin{enumerate}
\item For any objects $W_1$, $W_2$, $W_3$, $W_4$ and $M_1$ in ${\cal
C}$, any nonzero complex numbers $z_1$ and $z_2$ satisfying
$|z_1|>|z_2|>|z_1-z_2|>0$, any $P(z_1)$-intertwining map $I_1$ of type
${W_4}\choose {W_1 M_1}$ and $P(z_2)$-intertwining map $I_2$ of type
${M_1}\choose {W_2W_3}$ and any $w'_{(4)}\in W'_4$, we have that
$\gamma(I_1; 1_{W_1}, I_2)'(w'_{(4)})\in (W_1\otimes W_2\otimes
W_3)^{*}$ satisfies the $P^{(2)}(z_1-z_2)$-local grading restriction
condition.

\item For any objects $W_1$, $W_2$, $W_3$, $W_4$ and $M_2$ in ${\cal
C}$, any nonzero complex numbers $z_1$ and $z_2$ satisfying
$|z_1|>|z_2|>|z_1-z_2|>0$ and any $P(z_2)$-intertwining map $I^1$ of
type ${W_4}\choose {M_2 W_3}$ and $P(z_1-z_2)$-intertwining map $I^2$
of type ${M_2}\choose {W_1W_2}$ and any $w'_{(4)}\in W'_4$, we have
that $\gamma(I^1; I^2, 1_{W_3}) (w'_{(4)})\in (W_1\otimes W_2\otimes
W_3)^{*}$ satisfies the $P^{(1)}(z_2)$-local grading restriction
condition.
\end{enumerate}
\end{propo}

The proofs are analogous to those in \cite{tensor4}, since the
arguments continue to hold in the logarithmic case.

We will call either property in Proposition \ref{expansion} the {\it
expansion condition in the category ${\cal C}$}.

\begin{theo}\label{asso}
Assume that ${\cal C}$ is closed under $P(z)$-tensor products and that
the convergence and expansion conditions hold in ${\cal C}$. Then
there is a unique associativity isomorphism
\[
\alpha_{W_1 W_2 W_3}: (W_1\boxtimes_{P(z_1-z_2)}
W_2)\boxtimes_{P(z_2)} W_3 \longrightarrow W_1\boxtimes_{P(z_1)}
(W_2\boxtimes_{P(z_2)} W_3)
\]
for $W_1$, $W_2$ and $W_3$ in ${\cal C}$ such that
\begin{equation}\label{asso-elts}
\overline{\alpha}_{W_1 W_2 W_3}((w_{(1)}\boxtimes_{P(z_1-z_2)}
w_{(2)})\boxtimes_{P(z_2)} w_{(3)}) = w_{(1)}\boxtimes_{P(z_1)}
(w_{(2)}\boxtimes_{P(z_2)} w_{(3)})
\end{equation}
for all $w_{(1)}\in W_1$, $w_{(2)}\in W_2$ and $w_{(3)}\in W_3$,
where $\overline{\alpha}_{W_1 W_2 W_3}$ is the natural extension of
$\alpha_{W_1 W_2 W_3}$ to $\overline{(W_1\boxtimes_{P(z_1-z_2)}
W_2)\boxtimes_{P(z_2)} W_3}$.
\end{theo}
\noindent {\it Outline of proof}\hspace{2ex} Using the lemmas,
propositions and and theorems in this paper in place of those in
\cite{tensor4}, we can follow essentially the same steps as in
\cite{tensor4}, but treating the logarithmic subtleties with care, to
construct these associativity isomorphisms satisfying the condition
(\ref{asso-elts}).  \epfv

Theorem 14.8 in \cite{tensor4}, as generalized to our setting,
includes the statement that under the hypotheses of Theorem
\ref{asso}, the (general, nonmeromorphic) operator product expansion
in the generality of logarithmic intertwining operators exists, and
moreover, logarithmic intertwining operators satisfy the canonical
associativity property given in Proposition
\ref{associativityoflogintwops}.  If we drop the assumption that
${\cal C}$ is closed under $\boxtimes_{P(z)}$, then we still have the
existence of a (general, nonmeromorphic) operator product expansion in
the generality of logarithmic intertwining operators, in a suitable
sense.

\setcounter{equation}{0}
\setcounter{rema}{0}

\section{Differential equations and the convergence and extension
properties}

For the rest of this paper, for convenience we will take the grading
abelian groups $A$ and $\tilde A$ to be trivial, so that our (strongly
$A$-graded) conformal vertex algebra $V$ is simply a vertex operator
algebra (recall Remark \ref{trivialgroup}) and the (strongly graded)
generalized $V$-modules are the generalized $V$-modules in the sense
of Section 2 satisfying the grading restriction conditions referred to
in Definition \ref{def:dgv} and (\ref{set:dmltc}).  We also assume for
convenience that for any generalized $V$-module in ${\cal C}$, the
(generalized) weights form a discrete set of real numbers.

We first follow \cite{tensor4} to formulate the convergence and
extension properties in the sense of \cite{tensor4}, but now in the
logarithmic case, and to assert that these conditions together with
algebraic hypotheses imply the expansion condition.  We then state
results asserting that the theorems on differential equations in
\cite{diff-eqn} still hold.  This will imply that the condition called
$C_1$-cofiniteness in \cite{diff-eqn}, together with another algebraic
condition, implies that the convergence and extension properties
indeed hold.

Given objects $W_1$, $W_2$, $W_3$, $W_4$, $M_1$ and $M_2$ in ${\cal
C}$, let ${\cal Y}_1$, ${\cal Y}_2$, ${\cal Y}^1$ and ${\cal Y}^2$ be
intertwining operators of types ${W_4}\choose {W_1M_1}$, ${M_1}\choose
{W_2W_3}$, ${W_4}\choose {M_2W_3}$ and ${M_2}\choose {W_1W_2}$,
respectively. Consider the following conditions on the product of
${\cal Y}_1$ and ${\cal Y}_2$ and on the iterate of ${\cal Y}^1$ and
${\cal Y}^2$, respectively:

\begin{description}

\item[Convergence and extension property for products] There exists an
integer $N$ depending only on ${\cal Y}_1$ and ${\cal Y}_2$, and for
any $w_{(1)}\in W_1$, $w_{(2)}\in W_2$, $w_{(3)}\in W_3$, $w'_{(4)}\in
W'_4$, there exist $M\in\mathbb{N}$, $r_{k}, s_{k}\in \mathbb{R}$, $i_{k},
j_{k}\in \mathbb{N}$, $k=1,\dots,M$ and analytic functions
$f_{i_{k}j_{k}}(z)$ on $|z|<1$, $k=1, \dots, M$, satisfying
\[
\wt w_{(1)}+\wt w_{(2)}+s_{k}>N, \;\;\;k=1, \dots, M,
\]
such that
\[
\langle w'_{(4)}, {\cal Y}_1(w_{(1)}, x_2) {\cal Y}_2(w_{(2)},
x_2)w_{(3)}\rangle_{W_4} \lbar_{x_1= z_1, \;x_2=z_2}
\]
is convergent when $|z_1|>|z_2|>0$ and can be analytically extended to
the multi-valued analytic function
\[
\sum_{k=1}^{M}z_2^{r_{k}}(z_1-z_2)^{s_{k}}(\log z_2)^{i_{k}}
(\log(z_1-z_2))^{j_{k}}f_{i_{k}j_{k}}\left(\frac{z_1-z_2}{z_2}\right)
\]
in the region $|z_2|>|z_1-z_2|>0$.

\item[Convergence and extension property for iterates] There exists an
integer $\tilde{N}$ depending only on ${\cal Y}^1$ and ${\cal Y}^2$,
and for any $w_{(1)}\in W_1$, $w_{(2)}\in W_2$, $w_{(3)}\in W_3$,
$w'_{(4)}\in W'_4$, there exist $\tilde M\in\mathbb{N}$,
$\tilde{r}_{k}, \tilde{s}_{k}\in \mathbb{R}$, $\tilde{i}_{k},
\tilde{j}_{k}\in \mathbb{N}$, $k=1,\dots,\tilde M$ and analytic
functions $\tilde{f}_{\tilde{i}_{k}\tilde{j}_{k}}(z)$ on $|z|<1$,
$k=1, \dots, \tilde{M}$, satisfying
\[
\wt w_{(2)}+\wt w_{(3)}+\tilde{s}_{k}>\tilde{N}, \;\;\;k=1, \dots,
\tilde{M},
\]
 such that
\[
\langle w'_{(4)}, {\cal Y}^1({\cal Y}^2(w_{(1)}, x_0)w_{(2)},
x_2)w_{(3)}\rangle_{W_4} \lbar_{x_0=z_1-z_2,\;x_2=z_2}
\]
is convergent when $|z_2|>|z_1-z_2|>0$ and can be analytically
extended to the multi-valued analytic function
\[
\sum_{k=1}^{\tilde{M}} z_1^{\tilde{r}_{k}}z_2^{\tilde{s}_{k}} (\log
z_1)^{\tilde{i}_{k}} (\log
z)^{\tilde{j}_{k}}\tilde{f}_{\tilde{i}_{k}\tilde{j}_{k}}
\left(\frac{z_2}{z_1}\right)
\]
in the region $|z_1|>|z_2|>0$.
\end{description}

If for any objects $W_1$, $W_2$, $W_3$, $W_4$ and $M_1$ in ${\cal C}$
and any intertwining operators ${\cal Y}_1$ and ${\cal Y}_2$ of the
types as above, the convergence and extension property for products
holds, we say that {\it the convergence and extension property for
products holds in ${\cal C}$}.  We similarly define the meaning of the
phrase {\it the convergence and extension property for iterates holds
in ${\cal C}$}.

In the next theorem {\it only}, the term ``generalized module'' will
refer to this notion as introduced at the beginning of Section 2, with
no grading restriction conditions, that is, without the
strong-gradedness assumption (cf. the beginning of Section 4).

\begin{theo}\label{expansioncondnholds}
Suppose the following:
\begin{enumerate}
\item The category ${\cal C}$ is closed under $P(z)$-tensor products.
\item Every finitely-generated generalized $V$-module $W=\coprod_{n\in
\mathbb{R}}W_{[n]}$ that is lower truncated in the sense that
$W_{[n]}=0$ for $n$ sufficiently negative is an object in ${\cal C}$.
\item The convergence and extension property for either products or
iterates holds in ${\cal C}$.
\end{enumerate}
Then the expansion condition holds in ${\cal C}$ (recall Proposition
\ref{expansion}).
\end{theo}
\noindent {\it Outline of proof}\hspace{2ex}
This is proved by analogy with Theorem 16.2 in \cite{tensor4} and its
proof, in particular using Propositions
\ref{characterizationofbackslash}, \ref{5.12} and \ref{expansion}.
\epfv

\begin{rema}{\rm
As we commented at the end of Section 5, if the first of the three
assumptions in Theorem \ref{expansioncondnholds} is dropped, we still
have the existence of a (general, nonmeromorphic) operator product
expansion in the generality of logarithmic intertwining operators, in
a suitable sense.}
\end{rema}

Now we discuss generalizations of the results on differential
equations in \cite{diff-eqn}.

For our vertex operator algebra $V$, let $V_{+}=\coprod_{n>0}V_{(n)}$.
Let $W$ be a generalized $V$-module and let $C_{1}(W)= {\rm span
\,}\{u_{-1}w\;|\; u\in V_{+},\; w\in W\}$.  If $W/C_{1}(W)$ is finite
dimensional, we say that $W$ satisfies the {\it $C_{1}$-cofiniteness
condition}. If for any real number $N$, $\coprod_{n<N}W_{[n]}$ is
finite dimensional, we say that $W$ satisfies the {\it
quasi-finite-dimensionality condition}.

\begin{theo}\label{sys}
Let $n \geq 1$.  For $i=0, \dots, n+1$, let $W_{i}$ be a generalized
$V$-module satisfying the $C_{1}$-cofiniteness condition and the
quasi-finite-dimensionality condition.  Then for any $w_{(i)}\in
W_{i}$ ($i=0, \dots, n+1$), there exist $m\geq 1$ and
\[
a_{k, \;l}(z_{1}, \dots, z_{n})\in \mathbb{C}[z_{1}^{\pm 1}, \dots,
z_{n}^{\pm 1}, (z_{1}-z_{2})^{-1}, (z_{1}-z_{3})^{-1}, \dots,
(z_{n-1}-z_{n})^{-1}]
\]
($k=1, \dots, m$ and $l=1, \dots, n$) such that for any generalized
$V$-modules $\widetilde{W}_{i}$ ($i=1, \dots, n-1$) and any
logarithmic intertwining operators ${\cal Y}_{1}, {\cal Y}_{2}, \dots,
{\cal Y}_{n-1}, {\cal Y}_{n}$ of types ${W'_{0}\choose
W_{1}\widetilde{W}_{1}}$, ${\widetilde{W}_{1}\choose
W_{2}\widetilde{W}_{2}}, \dots, {\widetilde{W}_{n-2}\choose
W_{n-1}\widetilde{W}_{n-1}}$, ${\widetilde{W}_{n-1}\choose
W_{n}W_{n+1}}$, respectively, the (formal logarithmic) series
\begin{equation}
\langle w_{(0)}, {\cal Y}_{1}(w_{(1)}, z_{1})\cdots {\cal
Y}_{n}(w_{(n)}, z_{n})w_{(n+1)}\rangle
\end{equation}
satisfies the system of differential equations
\[
\frac{\partial^{m}\varphi}{\partial z_{l}^{m}}+ \sum_{k=1}^{m}a_{k,
\;l}(z_{1}, \dots, z_{n}) \frac{\partial^{m-k}\varphi}{\partial
z_{l}^{m-k}}=0,\;\;\;l=1, \dots, n
\]
and is absolutely convergent
in the region $|z_{1}|>\cdots >|z_{n}|>0$; such
$a_{k, \;l}(z_{1}, \dots, z_{n})$ can be chosen
so that the singular points of the system are
regular.
\end{theo}
\pf The result on the differential equations in \cite{diff-eqn} still
holds in our logarithmic case and the proof is essentially the same.
\epfv

We say that {\em products of arbitrary numbers of logarithmic
intertwining operators among objects in ${\cal C}$ are convergent and
extendable} if for any $n \geq 1$, any generalized $V$-modules $W_i$
for $i=0, \dots, n+1$ and $\widetilde{W}_i$ for $i=1, \dots, n-1$ in
${\cal C}$, and logarithmic intertwining operators ${\cal Y}_{1},
{\cal Y}_{2}, \dots, {\cal Y}_{n-1}, {\cal Y}_{n}$ of types
${W'_{0}\choose W_{1}\widetilde{W}_{1}}$, ${\widetilde{W}_{1}\choose
W_{2}\widetilde{W}_{2}}, \dots, {\widetilde{W}_{n-2}\choose
W_{n-1}\widetilde{W}_{n-1}}$, ${\widetilde{W}_{n-1}\choose
W_{n}W_{n+1}}$, respectively, each (formal logarithmic) series
\begin{equation} \langle w_{(0)}, {\cal Y}_{1}(w_{(1)}, z_{1})\cdots {\cal
Y}_{n}(w_{(n)}, z_{n})w_{(n+1)}\rangle
\end{equation}
is absolutely convergent in the region $|z_{1}|>\cdots> |z_{n}|>0$ and
can be analytically extended to a (multivalued) analytic function in
the region given by $z_{i}\ne 0$ for $i=1, \dots, n$ and $z_{i}\ne
z_{j}$ for $i, j=1, \dots, n$ with $i\ne j$.

We now have:
\begin{theo}\label{C_1pp}
Suppose that all generalized $V$-modules in ${\cal C}$ satisfy the
$C_{1}$-cofiniteness condition and the quasi-finite-dimensionality
condition. Then the convergence and extension properties for products
and iterates hold in ${\cal C}$ and in addition, products of arbitrary
numbers of logarithmic intertwining operators among objects in ${\cal
C}$ are convergent and extendable.
\end{theo}
\pf
This result follows immediately from Theorem \ref{sys}
and the theory of differential equations with
regular singular points.
\epfv

\setcounter{equation}{0}
\setcounter{rema}{0}

\section{Vertex and braided tensor category structures}

Finally, in this section we state the main theorems of our theory. For
the notion of vertex tensor category and the connection between these
categories and braided tensor categories, see \cite{tensorK}.

First we have:

\begin{theo}\label{conc}
Assume that the hypotheses of Theorem \ref{expansioncondnholds} hold
and in addition that products of arbitrary numbers of logarithmic
intertwining operators among objects in ${\cal C}$ are convergent and
extendable.  Then the category ${\cal C}$ has a natural structure of
vertex tensor category of central charge equal to the central charge
$c$ of $V$ such that for each $z\in \mathbb{C}^{\times}$, the tensor
product bifunctor $\boxtimes_{\psi(P(z))}$ associated to
$\psi(P(z))\in \tilde{K}^{c}(2)$ (see \cite{tensorK}) is equal to
$\boxtimes_{P(z)}$.  In particular, the category ${\cal C}$ has a
natural braided tensor category structure.
\end{theo}
\noindent {\it Outline of proof}\hspace{2ex}
By Remarks \ref{anyz_1} and \ref{bifunctor}, we have a tensor product
bifunctor for each nonzero complex number. General tensor product
bifunctors, associated to elements of the $\frac{c}{2}$-th power of
the determinant line bundle over the moduli space of spheres with one
negatively oriented puncture and two positively oriented punctures and
with local coordinates vanishing at the punctures, are constructed
using these bifunctors. The associativity isomorphisms for
$P(\cdot)$-tensor products have been constructed in Theorem
\ref{asso}. The commutativity isomorphisms for $P(\cdot)$-tensor
product functors are constructed using the map $\Omega_{-1}$ (recall
(\ref{Omega}) for the logarithmic intertwining operators associated to
$P(\cdot)$-tensor products. The general associativity and
commutativity isomorphisms are constructed from these special
isomorphisms.  The unit object is $V$ and the unit isomorphisms are
constructed using vertex operators defining the modules.  The
coherence properties follow easily from the convergence of the
products of logarithmic intertwining operators and the
characterizations of the associativity, commutativity and unit
isomorphisms in terms of tensor products of elements, such as
(\ref{asso-elts}).  \epfv

By Theorem \ref{C_1pp} we also have the following result, which yields
the vertex tensor category structure from purely algebraic hypotheses:
\begin{theo}
The conclusion of Theorem \ref{conc} is still true if the convergence
and extension property for products (or iterates) in ${\cal C}$, and
the condition that products of arbitrary numbers of logarithmic
intertwining operators among objects in ${\cal C}$ are convergent and
extendable, are replaced by the hypothesis that all generalized
$V$-modules in ${\cal C}$ are $C_{1}$-cofinite and quasi-finite
dimensional.
\end{theo}

Let $\kappa$ be a complex number not in $\mathbb{Q}_{\geq 0}$.
Consider the category ${\cal O}_{\kappa}$ of all modules for an affine
Lie algebra with central charge $\kappa-h$ (here $h$ is the dual
Coxeter number) of finite length, whose composition factors are simple
highest weight modules corresponding to weights which are dominant in
the direction of the finite-dimensional Lie algebra. Then the theorem
above and the results proved by the third author in \cite{Z} give the
following:

\begin{theo}
The category ${\cal O}_{\kappa}$ has a natural structure of vertex
tensor category, and in particular, braided tensor category.
\end{theo}

The braided tensor category structure in this theorem was first
constructed by Kazhdan and Lusztig (\cite{KL1}--\cite{KL5}) using a
very different method.

\bigskip

\noindent {\small \sc Department of Mathematics, Rutgers University,
Piscataway, NJ 08854}

\noindent {\em E-mail address}: yzhuang@math.rutgers.edu

\vspace{1em}

\noindent {\small \sc Department of Mathematics, Rutgers University,
Piscataway, NJ 08854}

\noindent {\em E-mail address}: lepowsky@math.rutgers.edu

\vspace{1em}

\noindent {\small \sc Department of Mathematics, Rutgers University,
Piscataway, NJ 08854}

\noindent {\em E-mail address}: linzhang@math.rutgers.edu

\end{document}